\documentclass{article}

\usepackage[a4paper]{geometry}
\usepackage[utf8]{inputenc}
\usepackage[T1]{fontenc}
\usepackage{ae}
\usepackage{eucal}
\renewcommand{\geq}{\geqslant}
\renewcommand{\leq}{\leqslant}
\renewcommand{\phi}{\varphi}
\usepackage{amsmath,amsfonts,amssymb}
\usepackage{hyperref}
\usepackage{url}
\usepackage{theorem}
\newtheorem{theorem}{Theorem}
\newtheorem{lemma}[theorem]{Lemma}
\newtheorem{definition}[theorem]{Definition}
\newtheorem{remark}[theorem]{Remark}
{
\theoremstyle {break}
\newtheorem {algorithm}[theorem]{Algorithm}
}
\newcommand {\inoutput}[2]{\textsc {Input:} #1\\
   \textsc {Output:} #2}
\newcommand {\keyword}[1]{\textbf {#1}}

\newcommand {\weil}{e_n}
\newcommand {\rweil}{e_r}
\newcommand {\tate}{e_n^\mathrm {T}}
\newcommand {\rtate}{e_r^\mathrm {T}}
\newcommand {\tatered}{e_n^\mathrm {T'}}

\newcommand {\ate}{e_r^\mathrm {A}}
\newcommand {\twistedate}{e_r^\mathrm {\tilde A}}
\newcommand {\atei}{e_r^{\mathrm {A}_i}}
\newcommand {\twistedatei}{e_r^{\mathrm {\tilde A}_i}}
\newcommand {\hess}{e_r^\mathrm {H}}
\newcommand {\twistedhess}{e_r^\mathrm {\tilde H}}
\newcommand {\fhess}{f_{t, y, Q}}
\newcommand {\ftwistedhess}{f_{t, y, P}}

\newcommand {\Z}{\mathbb Z}
\newcommand {\F}{\mathbb F}

\newcommand {\Oc}{\mathcal O}
\newcommand {\id}{\mathrm {id}}
\newcommand {\ord}{\operatorname {ord}}
\newcommand {\supp}{\operatorname {supp}}
\newcommand {\ddiv}{\operatorname {div}}
\newcommand {\kernel}{\operatorname {ker}}
\newcommand {\Norm}{\operatorname {N}}
\newcommand {\Tr}{\operatorname {Tr}}
\newcommand {\Div}{\operatorname {Div}}
\newcommand {\Prin}{\operatorname {Prin}}
\newcommand {\Pic}{\operatorname {Pic}}
\newcommand {\lc}{\operatorname {lc}}
\newcommand {\Aut}{\operatorname {Aut}}
\newcommand {\charac}{\operatorname {char}}

\newenvironment {proofof}[1]
   {{\parindent0cm \bf Proof#1:}}
   {\hspace* {\fill} $\Box$}
\newenvironment {proof}
   {\begin {proofof}{}}
   {\end {proofof}}

\title {Bilinear pairings on elliptic curves}
\author {Andreas Enge\footnote {INRIA, LFANT, F-33400 Talence, France \newline
CNRS, IMB, UMR 5251, F-33400 Talence, France \newline
Univ. Bordeaux, IMB, UMR 5251, F-33400 Talence, France}}
\date {14 February 2014}

\begin {document}

\maketitle

\begin {abstract}
We give an elementary and self-contained introduction to pairings on
elliptic curves over finite fields.
The three different definitions of the Weil pairing that can be found in
the literature are stated and proved to be equivalent using Weil reciprocity.
Pairings with shorter loops, such as the ate, ate$_i$, R-ate and optimal
pairings, together with their twisted variants, are presented with proofs
of their bilinearity and non-degeneracy.
Finally, we review different types of pairings in a cryptographic context.
This article can be seen as an update chapter to
A. Enge,
\textit {Elliptic Curves and Their Applications to Cryptography --
An Introduction}, Kluwer Academic Publishers 1999.
\end {abstract}

\section {Introduction}
\label {sec:intro}

Consider three abelian groups $G_1$, $G_2$ (written additively) and $G_3$
(written multiplicatively), which can equivalently
be seen as $\Z$-modules. A \textit {pairing} on $G_1$ and $G_2$ with values
in $G_3$ is a $\Z$-bilinear map
\[
e : G_1 \times G_2 \to G_3,
\]
so that
\[
e (aP, bQ) = e (P, Q)^{ab}
\]
for all elements $P \in G_1$, $Q \in G_2$ and integers $a$ and $b$.
In the following, $G_1$ and $G_2$ will be groups related to an elliptic
curve $E$ defined over some field $K$: They will be subgroups of the
elliptic curve group (in the case of the Weil pairing of \S\ref {sec:weil})
or subgroups and quotient groups (in the case of the Tate pairing of
\S\ref {sec:tate} and related pairings presented in \S\ref {sec:shortloop}).
The group $G_3$ will be a subgroup or a quotient of the multiplicative
group $K^*$.

Elliptic curve cryptosystems are currently among the most efficient
public-key systems. Their security relies on the difficulty of computing
discrete logarithms in suitable instances of elliptic curves over
finite fields, that is, on the difficulty of computing $x$ given two
points $P$ and $R = xP$ on the curve.
Pairings then transport the discrete logarithm problem from the curve
into the multiplicative group of a finite field, where it is potentially
easier to solve \cite {Odlyzko13}:
As $e (R, Q) = e (P, Q)^x$, the discrete logarithm
of $e (R, Q)$ with respect to the basis $e (P, Q)$ yields~$x$.
Consequently, pairings have first been suggested as a means
of attacking elliptic curve cryptosystems
\cite{MeOkVa93,FrRu94}. First constructive cryptographic applications have
been described in \cite{Joux00,SaOhKa00,BoFr01}, and since then, the number
of publications introducing pairing-based cryptographic primitives has
exploded. A new conference series, \textit {Pairing}, is devoted to
the topic \cite{TaOkOkOk07,GaPa08,ShWa09,JoMiOt10,AbLa13,CaZh14}.

This document provides a self-contained introduction to pairings and aims at
summarising the state of the art as far as the definitions of different
pairings and their cryptographic use are concerned. While being as accessible
as possible, we do not sacrifice mathematical rigour, in the style of
\cite {Enge99}, of which the current article can be seen as an update
chapter.
While most of the following holds over arbitrary perfect or even more general
fields, we limit the presentation to
the only case of interest in the cryptographic context, namely $K$ being a
finite field $\F_q$ with $q$ elements. Pairings can be defined in Jacobians
of arbitrary curves or, more generally, in abelian varieties. However, due to
recent progress in solving the discrete logarithm problem (see the survey
\cite{Enge08}), only elliptic curves and hyperelliptic curves of genus~$2$
appear to be suited for cryptography. For the latter, the problem of
finding instances in which the pairing has values in a sufficiently
small finite field to be efficiently computed (see the definition of the
embedding degree at the beginning of~\S\ref {sec:weil}) and in which
the size of the subgroup is relatively close to that of the full group to
allow for bandwidth-efficient protocols has not yet been solved in a
satisfactory way. So in the following we consider only elliptic curves.

An excellent survey is given by Galbraith in \cite{Galbraith05}. We
complement his presentation by concentrating on the Weil pairing instead of
the Tate pairing and by reporting on progress made after the publication of
\cite {Galbraith05} concerning pairings with shorter evaluation loops.

\section {Elliptic curves and Weil reciprocity}
\label {sec:divisors}

\subsection {Divisors and group law}
\label {ssec:divisors}

We assume the reader to be familiar with basic algebra, in particular with
finite fields. For proofs of the following facts on elliptic curves, see
\cite{Silverman86,Enge99}. Other sources for the use of elliptic curves in
cryptography are \cite{CoFrAvDoLaNgVe06,BlSeSm99}.
From now on, we assume that $K = \F_q = \F_{p^m}$ is the finite field of
characteristic $p$ with $q$ elements.
(This is motivated by the cryptologic applications and meant to ease
the exposition. All statements concerning the Weil pairing hold in fact over
arbitrary fields. The definition given of the Tate pairing
in~\S\ref {sec:tate}, however, is not valid for all fields;
over finite fields, it yields a non-degenerate pairing.)

In several places, we will consider the algebraic closure $\overline K$
for convenience; this could be replaced by a sufficiently large extension
field to contain the coordinates of all points under consideration. An
\textit {elliptic curve} over $K$ is given by a non-singular, absolutely
irreducible \textit {long Weierstra{\ss} equation}
\[
E : Y^2 + (a_1 X + a_3) Y = X^3 + a_2 X^2 + a_4 X + a_6
\]
with $a_i \in K$. If $p \geq 5$, the equation can be transformed into short
Weierstra{\ss} form in which all but $a_4$ and $a_6$ vanish. The points on
$E$ are given by the \textit {affine} points $(x, y) \in K^2$ satisfying the
equation, together with a \textit {projective point at infinity} $\Oc$. The
\textit {coordinate ring} of $E$ is the ring $K [E] = K [X, Y] / (E)$ of
polynomial functions, its \textit {function field} $K (E) = K (X)[Y] / (E) =
\{ a (X) + b (X) Y : a, b \in K (X) \}$ is the set of \textit {rational
functions} from $E$ to $K \cup \{ \infty \}$; the value $\infty$ is reached
when the function has a pole in a point. It turns out that the points on $E$
are in a one-to-one correspondence with the discrete valuation rings of $K
(E)$, given by the rings $\Oc_P$ of functions that do not have a pole in~$P$.

The set $E (K)$ of points on $E$ with coordinates in $K$ (including $\Oc$)
can be turned into a finite abelian group via the tangent-and-chord law:
$\Oc$ is the neutral element of the group law, and three points on a line sum
to $\Oc$. The only delicate point in proving the group law is associativity;
the simplest proof, which also generalises to other curves, is sketched in
the following. It uses divisors,
which are needed anyway to define pairings. So let
\[
\Div (E) = \left\{\sum_P n_P [P] : P \in E (K), n_P \in \Z,
\text {only finitely many $n_P$ are non-zero} \right\}
\]
be the free abelian group over the points on $E$, define the degree of a
divisor as the sum $\sum n_P$ of its coefficients, and let $\Div^0 (E)$ be
the subgroup of $\Div (E)$ consisting of divisors of degree~$0$. To a
rational function $f \in K (E)$, associate its divisor $\ddiv (f) = \sum_P
\ord_P (f) [P]$, where $\ord_P (f)$ is the valuation of $f$ with respect to
$\Oc_P$, that is, $\ord_P (f) > 0$ if $P$ is a zero of $f$, $\ord_P (f) < 0$
if $P$ is a pole, and $\ord_P (f) = 0$ otherwise. Let $\Prin (E) = \{ \ddiv
(f) : f \in K (E) \} \subseteq \Div^0 (E)$ be the set of \textit {principal
divisors}. Then the quotient $\Pic^0 (E) = \Div^0 (E) / \Prin (E)$ is
evidently a group, and it can be identified with $E (K)$ via $P \mapsto [P] -
[\Oc]$, which maps $\Oc$ to the neutral element~$O$.

Let $\sim$ denote equivalence modulo $\Prin^0 (E)$. The geometric
tangent-and-chord law is recovered as follows. For a point $R = (x_R, y_R)$,
let
\begin {equation}
\label {eq:v}
v_R = X - x_R
\end {equation}
be the vertical line through $R$. Then $\ddiv (v_R) = [R] + [\overline R] - 2
[\Oc] \sim 0$ with $\overline R = (x_R, -y_R-a_1 x_R - a_3)$,
so that $-R = \overline R$. For
two points $P = (x_P, y_P)$ and $Q = (x_Q, y_Q)$ with $Q \neq -P$ let
$\ell_{P, Q}$ be the chord through $P$ and $Q$ if $P \neq Q$ or the tangent
at $P$ if $P = Q$:
\begin {equation}
\label {eq:l}

\begin {array}{lll}
\lambda_{P, Q} & = & \left\{
\begin {array}{ll}
\frac {y_Q - y_P}{x_Q - x_P}
& \text { if } P \neq Q \\
\frac {3 x_P^2 + 2 a_2 x_P + a_4}{2 y_P + a_1 x_P + a_3}
& \text { if } P = Q
\end {array}
\right. \\
\ell_{P, Q} & = & (Y - y_P) - \lambda_{P, Q} (X - x_P)
\end {array}
\end {equation}
Then $\ell_{P, Q}$ intersects $E$ in a third point $R = (x_R, y_R) \neq \Oc$,
and $\ddiv \left( \frac {\ell_{P, Q}}{v_R} \right) = \ddiv (\ell_{P, Q}) -
\ddiv (v_R) = \big( [P] + [Q] + [R] - 3 [\Oc] \big) - \big( [R] + [\overline
R] - 2 [\Oc] \big) = [P] + [Q] - [\overline R] - [\Oc] \sim 0$ implies that
$P + Q = \overline R$.

By induction, this proves the following characterisation of principal
divisors.

\begin {theorem}
\label {th:principality}
A divisor $D = \sum_P n_P [P]$ is principal if and only if $\deg D = 0$ and
$\sum_P n_P P = \Oc$ on $E$. The function associated to a principal divisor
is unique up to multiplication by constants in $K^\ast$.
\end {theorem}

It is often useful to assume the following normalisation.

\begin {definition}
\label {def:monic}
The \textit {leading coefficient} of a function $f$ at $\Oc$ is
\[
\lc (f) =
\left(\left( \frac {X}{Y} \right)^{-\ord_\Oc (f)}f \right) (\Oc).
\]
A function $f$ is \textit {monic} at $\Oc$ if $\lc (f) = 1$.
\end {definition}
In particular, the lines $v_R$ and $\ell_{P, Q}$ given above for the
tangent-and-chord law are monic at~$\Oc$, and this implies that the
functions computed in Algorithm~\ref {alg:miller} will also be
monic at~$\Oc$.

\subsection {Rational maps, isogenies and star equations}

Let $E$, $E'$ be two elliptic curves over the same field $K$. A \textit
{rational map} $\alpha:E\to E'$ is an element of $E'(K (E))$. Explicitly,
$\alpha$ is given by rational functions in $X$ and $Y$ that satisfy the
Weierstra{\ss} equation for $E'$. Unless $\alpha$ is constant, it is
surjective. If $\alpha (\Oc) = \Oc'$, then $\alpha$ is in fact a group
homomorphism, and it is called an \textit {isogeny}. If furthermore $E = E'$,
then $\alpha$ is called an \textit {endomorphism}. The endomorphisms that are
most important in the following are multiplications by an integer $n$,
denoted by $[n]$.

A non-constant rational map $\alpha : E \to E'$ induces an injective
homomorphism of function fields $\alpha^\ast : K (E') \to K (E)$, $f' \mapsto
f' \circ \alpha$; the \textit {degree} of $\alpha$ is the degree of the
function field extension $[K (E) : \alpha^\ast (K (E'))]$. For instance,
$\deg ([n]) = n^2$. If $\alpha$ is an isogeny, there is another isogeny $\hat
\alpha$ of the same degree, called its \textit {dual}, such that $\hat \alpha
\circ \alpha = [\deg \alpha]$.

For a point $P \in E$ and $P' = \alpha (P)$, there is an integer $e_\alpha
(P)$, called \textit {ramification} index, such that $\ord_P (\alpha^\ast
(f')) = e_\alpha (P) \ord_{P'} (f')$ for any $f' \in K (E')$. When $\alpha$
is an isogeny, $e_\alpha (P)$ is independent of $P$. In this case, we have
$\deg \alpha = e_\alpha \cdot \#(\ker \alpha)$, and two extreme cases can
occur: If $e_\alpha = 1$, then $\alpha$ is called \textit {separable}; in
particular, $[n]$ is separable if $p \nmid n$. If $\#(\ker \alpha) = 1$, then
$\alpha$ is (up to isomorphisms) a power of the \textit {purely inseparable}
\textit {Frobenius endomorphism} $(x, y) \mapsto (x^q, y^q)$ of degree and
ramification index~$q$. An arbitrary isogeny can be decomposed into a
separable one and a power of Frobenius, which is often convenient for proving
theorems.

The ramification index allows to define a homomorphism $\alpha^\ast : \Div
(E') \to \Div (E)$ on divisors by
\[
\alpha^\ast ([P']) = \sum_{P \in \alpha^{-1} (P')} e_\alpha (P) [P]
\]
in such a way that the maps $\alpha^\ast$ on functions and divisors are
compatible; the proof follows immediately from the definition of $e_\alpha$.

\begin {theorem}[Upper star equation]
\label {th:upperstar}
If $\alpha : E \to E'$ is a non-constant rational map and $f' \in K (E')$,
then
\[
\alpha^\ast (\ddiv (f')) = \ddiv (\alpha^\ast (f')).
\]
\end {theorem}

On the other hand, the map $\alpha_\ast : \Div (E) \to \Div (E')$ is defined
by $\alpha_\ast ([P]) = [\alpha (P)]$. A corresponding map on function fields
$K (E) \to K (E')$ can be defined by
\[
\alpha_\ast (f) = (\alpha^\ast)^{-1} \left( \Norm_{K (E) / \alpha^\ast (K
(E'))} (f) \right),
\]
where $\Norm$ denotes the norm with respect to the function field extension.
The map $\alpha_\ast$ is well-defined because the norm is an element of
$\alpha^\ast (K (E'))$, so that a preimage exists, and because $\alpha^\ast$
is injective, so that the preimage is unique.

It is shown in \cite[(18)]{ChCo90} that
\begin {equation}
\label {eq:norm}
\Norm_{K (E) / \alpha^\ast (K (E'))} (f)
= \left( \prod_{R \in \ker \alpha} (f \circ \tau_R) \right)^{e_\alpha},
\end {equation}
where $\tau_R$ is the translation by $R$; the product accounts for the
separable, the exponent for the inseparable part of the isogeny. This can be
used to show the following result:

\begin {theorem}[Lower star equation]
\label {th:lowerstar}
If $\alpha : E \to E'$ is a non-constant rational map and $f \in K (E)$, then
\[
\alpha_\ast (\ddiv (f)) = \ddiv (\alpha_\ast (f)).
\]
\end {theorem}

\subsection {Weil reciprocity}

The key to the definition of pairings is the evaluation of rational functions
in divisors. For $D = \sum_P n_P [P]$ let its \textit {support} be $\supp (D)
= \{ P : n_P \neq 0 \}$. The evaluation of a rational function $f$ in points
is extended to a group homomorphism from divisors (with support disjoint from
$\supp (\ddiv f)$) to $K^\ast$ via
\[
f \left( \sum_P n_P [P] \right)
= \prod_P f(P)^{n_P}.
\]

In order to handle common points in the supports, let the \textit {tame
symbol} of two functions $f$ and $g \in K (E)$ be defined as
\[
\langle f, g \rangle_P = (-1)^{\ord_P (f) \ord_P (g)}
\left( \frac {f^{\ord_P (g)}}{g^{\ord_P (f)}} \right) (P).
\]
\begin {theorem}[Generalised Weil reciprocity]
\label {th:reciprocity}
If $f$, $g \in K (E)$, then
\[
\prod_{P \in E (\overline K)} \langle f , g \rangle_P = 1.
\]
In particular, if $\supp (f) \cap \supp (g) = \emptyset$, then
\begin {equation}
\label {eq:weilrecdisjoint}
f (\ddiv g) = g (\ddiv f).
\end {equation}
\end {theorem}
For a proof, see \cite[\S7]{ChCo90}.

\section {Weil pairing}
\label {sec:weil}

Let $E [n] = \{ P \in E (\overline K) : n P = \Oc \} = \kernel ([n])$ be the
set of $n$-torsion points of $E$, which are in general not defined over $K$
itself.
For future reference, we denote by $E (K)[n] = E [n] \cap E (K)$ the set of
points of $E [n]$ defined over $K$, which contains at least~$\Oc$.
From now on, we will assume that
$\gcd (n, p) = 1$; then the group $E [n]$ is finite
and isomorphic to $\Z / n \Z \times \Z / n \Z$. The field $L$ obtained by
adjoining to $K = \F_q$ all coordinates of $n$-torsion points is thus a
finite field extension $\F_{q^k}$, and $k$ is called the \textit {embedding
degree} of the $n$-torsion and $\F_{q^k}$ its \textit {embedding field}.
We have $L \supseteq K (\zeta_n)$, where
$\zeta_n$ is a primitive $n$-th root of unity, and equality holds in the case
of main cryptographic interest, namely that $n$ is a prime and $n \nmid q -
1$ by \cite[Th.~1]{BaKo98}. Then $k$ is the smallest integer such that $n
\mid q^k - 1$.

\begin {theorem}
\label {th:weil}
The Weil pairing is a map
\[
\weil : E [n] \times E [n] \to \mu \subset L^\ast,
\]
where $\mu$ is the set of $n$-th roots of unity in $L$, satisfying the
following properties:
\begin {enumerate}
\item
Bilinearity:
\begin {eqnarray*}
&& \weil (P_1 + P_2, Q) = \weil (P_1, Q) \weil (P_2, Q), \\
&& \weil (P, Q_1 + Q_2) = \weil (P, Q_1) \weil (P, Q_2)
\quad \forall P, P_1, P_2, Q, Q_1, Q_2 \in E [n];
\end {eqnarray*}
\item
Identity:
\[
\weil (P, P) = 1 \quad \forall P \in E [n];
\]
\item
Alternation:
\[
\weil (P, Q) = \weil (Q, P)^{-1} \quad \forall P, Q \in E [n];
\]
\item
Non-degeneracy:
For any $P \in E [n] \backslash \{ \Oc \}$, there is a $Q \in E [n]$, and
for any $Q \in E [n] \backslash \{ \Oc \}$, there is a $P \in E [n]$
such that $\weil (P, Q) \neq 1$;
\item
Compatibility with isogenies:
\begin {eqnarray*}
\weil (\alpha (P), \alpha (Q)) & = & \weil (P, Q)^{\deg \alpha}, \\
\weil (P', \alpha (Q))         & = & \weil (\hat \alpha (P'), Q)
\end {eqnarray*}
for $P$, $Q \in E [n]$, $P' \in E' [n]$ and $\alpha : E \to E'$ a non-zero
isogeny defined over $L$.
In particular, $\alpha$ may be the Frobenius endomorphism on~$E$
of degree~$q$.
\end {enumerate}
\end {theorem}

In the literature, there are in fact three equivalent definitions of the Weil
pairing, and depending on which one is chosen, the different properties are
more or less easy to prove, the most intricate one being non-degeneracy.
In the following, we show equivalence of these definitions, which is
also non-trivial and makes intensive use of Weil reciprocity, and we
prove the five properties of the Weil pairing using for each the definition
that yields the easiest proof.

\paragraph {First definition of the Weil pairing
(\cite[{\S}III.8]{Silverman86},\cite[\S3.7]{Enge99}).}
For $P \in E [n]$, consider $D = [n]^\ast ([P] - [\Oc]) = \sum_{R \in E [n]}
([P_0 + R] - [R])$, where $P_0$ is any point such that $n P_0 = P$. By
Theorem~\ref {th:principality}, $D$ is principal; let $g_P$ be such that
$\ddiv g_P = D$. Let again $\tau_Q : R \mapsto R + Q$ denote the translation
by $Q \in E [n]$. Then
\begin {equation}
\label {eq:weil1}
\weil (P, Q) = \frac {g_P \circ \tau_Q}{g_P}.
\end {equation}
While $g_P$ is defined only up to multiplication by non-zero constants, the
quotient is a well-defined rational function. Since $\ddiv (g_P \circ \tau_Q)
= \ddiv (\tau_Q^\ast (g_P)) = \tau_Q^\ast (\ddiv g_P)$ by Theorem~\ref
{th:upperstar} and the latter divisor equals
\[
\sum_{R \in E [n]} \left( [P_0 + R - Q] - [R - Q] \right) = \ddiv g_P
\]
for $Q \in E [n]$, the Weil pairing yields indeed a constant in $\overline
K$. That it yields an $n$-th root of unity follows from bilinearity

\begin {proofof}{ of Theorem~\ref {th:weil}(a)}
Using (c), proved below, it is sufficient to show linearity in the second
argument, which follows from the definition:
\begin {eqnarray*}
\weil (P, Q_1 + Q_2)
& = & \frac {g_P \circ \tau_{Q_1 + Q_2}}{g_P}
= \left( \frac {g_P \circ \tau_{Q_1}}{g_P} \circ \tau_{Q_2} \right)
\frac {g_P \circ \tau_{Q_2}}{g_P} \\
& = & \weil (P, Q_1) \weil (P, Q_2)
\text { since the constant } \weil (P, Q_1) \\
&& \text {is invariant under } \tau_{Q_2}.
\end {eqnarray*}
\end {proofof}

\begin {proofof}{ of Theorem~\ref {th:weil}(d)}
We sketch the approach of \cite[Prop.~3.60]{Enge99}. Using (c), it is
sufficient to show non-degeneracy with respect to the first argument. For $P
\in E [n]$, suppose that $\weil (P, Q) = 1$ for all $Q \in E [n]$. This means
that $g_P$ is invariant under translations by all $Q \in E [n] = \kernel
([n])$, so that all conjugates of $g_P$ with respect to the field extension
$K (E) / [n]^\ast (K (E))$ are $g_P$ itself, see \eqref {eq:norm}. Hence,
there is a function $f_P$ such that $g_P = [n]^\ast (f_P)$. By Theorem~\ref
{th:upperstar}, this implies that $\ddiv f_P = [P] - [\Oc]$, which by
Theorem~\ref {th:principality} implies $P = \Oc$.
\end {proofof}

\begin {proofof}{ of Theorem~\ref {th:weil}(e)}
As a homomorphism, $\alpha$ commutes with $[n]$, and being surjective,
it acts as a permutation on $E [n]$. So
\begin {eqnarray*}
\ddiv (g_{\alpha (P)})
& = & \sum_{R \in E [n]} \left( [\alpha (P_0) + R] - [R] \right) \\
& = & \sum_{S \in E [n]} \left( [\alpha (P_0) + \alpha (S)]
- [\alpha (S)] \right)
\text { where } R = \alpha (S) \\
& = & \alpha_\ast (\ddiv (g_P)) \\
& = & \ddiv (\alpha_\ast (g_P))
\text { by Theorem~\ref {th:lowerstar}}.
\end {eqnarray*}
This implies
$g_{\alpha (P)} = c \alpha_\ast (g_P)$ for some $c \in K^\ast$,
and
\begin {equation}
\label {eq:divisorisogeny}
g_{\alpha (P)} \circ \alpha = \alpha^\ast (g_{\alpha (P)})
= c \left( \prod_{R \in \ker \alpha} (g_P \circ \tau_R) \right)^{e_\alpha}
\end {equation}
by \eqref {eq:norm}.
Hence,
\begin {eqnarray*}
\weil (\alpha (P), \alpha (Q))
& = & \weil (\alpha (P), \alpha (Q)) \circ \alpha
= \frac {g_{\alpha (P)} \circ \tau_{\alpha (Q)}}{g_{\alpha (P)}}
  \circ \alpha
= \frac {g_{\alpha (P)} \circ \alpha \circ \tau_Q}
  {g_{\alpha (P)} \circ \alpha} \\
& = & \left( \prod_{R \in \ker \alpha}
\left( \frac {g_P \circ \tau_Q}{g_P} \right)
\circ \tau_R \right)^{e_\alpha}
= \weil (P, Q)^{e_\alpha \cdot \#(\ker \alpha)} \\
& = & \weil (P, Q)^{\deg \alpha}.
\end {eqnarray*}
Concerning the second equation, let $P$ be such that $\alpha (P) = P'$; then
$\hat \alpha (P') = (\hat \alpha \circ \alpha) (P) = (\deg \alpha) P$, and
\[
\weil (\hat \alpha (P'), Q)
= \weil (P, Q)^{\deg \alpha}
= \weil (\alpha (P), \alpha (Q))
= \weil (P', \alpha (Q)).
\]
\end {proofof}

\paragraph {Second definition of the Weil pairing.}
For $P, Q \in E [n] \backslash \{ \Oc \}$, $P \neq Q$,
let $f_P$ and $f_Q$ be such that
$\ddiv f_P = n [P] - n [\Oc]$ and
$\ddiv f_Q = n [Q] - n [\Oc]$, which is possible by
Theorem~\ref {th:principality}. Then
\begin {equation}
\label {eq:weil2}
\weil (P, Q) = (-1)^n \cdot \frac {f_P (Q)}{f_Q (P)}
   \cdot\frac {f_Q}{f_P} (\Oc);
\end {equation}
if $f_P$ and $f_Q$ are chosen monic at $\Oc$ as in
Definition~\ref {def:monic}, then
\[
\weil (P, Q) = (-1)^n \cdot \frac {f_P (Q)}{f_Q (P)}.
\]
For $P = Q$ or one or both of $P$ and $Q$ being $\Oc$, the definition
needs to be completed by $\weil (P, Q) = 1$.

\begin {remark}
\label {rem:weil2}
This definition is the most suited one for computations, see
Algorithm~\ref {alg:miller}.
The factor $(-1)^n$ is often missing in the literature.
\end {remark}

\begin {proofof}{ of equivalence of the two definitions}
We essentially follow \cite[\S10]{ChCo90}.
Assume that $\weil$ is defined as in
\eqref {eq:weil1}.

Let $P_0$ and $Q_0$ be such that $n P_0 = P$ and $n Q_0 = Q$.
Let $g_P$ be the function, monic at $\Oc$, such that
\[
\ddiv (g_P) = \sum_{R \in E [n]} \left( [P_0 + R] - [R] \right),
\]
and similarly for $g_Q$.

If $P = \Oc$, we may take $P_0 = \Oc$, which shows that $g_\Oc = 1$
and $\weil (\Oc, Q) = 1$. If $Q = \Oc$, then $\tau_Q = \id$, and
$\weil (P, \Oc) = 1$. So from now on, $P$, $Q \neq \Oc$.

Let $h_Q$ be the function, monic at $\Oc$, such that
\[
\ddiv h_Q = (n - 1) [Q_0] + [Q_0 - Q] - n [\Oc],
\]
which exists by Theorem~\ref {th:principality},
and let $H_Q = \prod_{R \in E [n]} (h_Q \circ \tau_R)$.
By comparing divisors and leading coefficients,
$H_Q = \lc (H_Q) \cdot g_Q^n$.

By generalised Weil reciprocity of Theorem~\ref {th:reciprocity},
we have
\[
\prod_{S \in \supp (\ddiv g_P) \cup \supp (\ddiv h_Q)}
\langle g_P, h_Q \rangle_S = 1.
\]

If $P \neq Q$, then
$\supp (\ddiv g_P) \cap \supp (\ddiv h_Q) = \{ \Oc \}$,
and we easily compute the different contributions of tame symbols:
\begin {eqnarray*}
\langle g_P, h_Q \rangle_{Q_0}
& = & g_P^{n-1} (Q_0) \\
\langle g_P, h_Q \rangle_{Q_0 - Q}
& = & g_P (Q_0 - Q) \\
\langle g_P, h_Q \rangle_{P_0 + R}
& = & h_Q^{-1} (P_0 + R) \text { for } R \in E [n] \\
\langle g_P, h_Q \rangle_R
& = & h_Q (R) \text { for } R \in E [n] \backslash \{ \Oc \} \\
\langle g_P, h_Q \rangle_\Oc
& = & (-1)^n \frac {h_Q}{g_P^n} (\Oc) = (-1)^n
\text { since $g_P$ and $h_Q$ are monic at $\Oc$}.
\end {eqnarray*}
Multiplying them together, we find that
\begin {eqnarray*}
1 & = & g_P^n (Q_0)
\underbrace {\frac {g_P (Q_0 - Q)}{g_P (Q_0)}}_{\:\:
\frac {g_P}{g_P \circ \tau_Q} (Q_0 - Q) = \weil (P, Q)^{-1} \:\:}
\underbrace {\frac {1}{H_Q (P_0)}}_{\lc (H_Q)^{-1} g_Q (P_0)^{-n}}
\underbrace {\frac {H_Q}{h_Q} (\Oc)}_{\lc (H_Q)}
(-1)^n \\
& = &
(-1)^n \, \frac {g_P^n (Q_0)}{g_Q^n (P_0)} \cdot\frac {1}{\weil (P, Q)}.
\end {eqnarray*}
Since $\ddiv (g_P^n) = n [n]^\ast ([P] - [\Oc]) = [n]^\ast \ddiv (f_P)$,
Theorem~\ref {th:upperstar} implies that
\[
g_P^n = c^{-1} \cdot [n]^\ast (f_P)
\]
with $c = \lc ([n]^\ast (f_P)) =
\left( (f_P \circ [n]) \frac {X^n}{Y^n} \right) (\Oc)$.
An analogous equation holds for $g_Q^n$,
so that
\[
\frac {g_P^n (Q_0)}{g_Q^n (P_0)}
= \frac {f_P (Q)}{f_Q (P)} \cdot \frac {f_Q}{f_P} (\Oc).
\]

If $P = Q$, then $\supp (\ddiv (h_Q)) \subseteq \supp (\ddiv (g_Q))$,
and a similar computation shows that $\weil (P, P) = 1$.
\end {proofof}

\begin {proofof}{ of Theorem~\ref {th:weil}(b)}
This is part of the second definition. (The only statement needing proof
is that this also holds for the first definition, as shown above.)
\end {proofof}

\begin {proofof}{ of Theorem~\ref {th:weil}(c)}
This is immediate from \eqref {eq:weil2}.
\end {proofof}

\paragraph {Third definition of the Weil pairing.}
For any degree zero divisor $D$ such that $n D \sim 0$ in $\Pic^0 (E)$, we
denote by $f_D$ the function. monic at~$\Oc$, such that
$\ddiv (f_D) = n D$; thus $f_{[P] - [\Oc]} = f_P$.
Choose $D_P \sim [P] - [\Oc]$ and $D_Q \sim [Q] - [\Oc]$ with disjoint
supports. Then
\begin {equation}
\label {eq:weil3}
\weil (P, Q) = \frac {f_{D_P} (D_Q)}{f_{D_Q} (D_P)}.
\end {equation}
Note the similarity with \eqref {eq:weil2}, but also the missing factor
$(-1)^n$, due to the common pole~$\Oc$ of $f_P$ and $f_Q$.

\begin {remark}
The third definition corresponds to Weil's original one in
\cite {Weil40}. The first definition is given in \cite{Silverman86,Enge99}
with the roles of $P$ and $Q$ exchanged, which by the alternation property
yields the inverse of the Weil pairing. The definition with $P$ and $Q$
in the order of this paper is used in the Notes on Exercises, p.~462 of the
second edition of \cite {Silverman86}, as well as in~\cite {Silverman13}.
\end {remark}

One needs to check that \eqref {eq:weil3} is well-defined.
Let $D_Q' \sim [Q] - [\Oc]$ be another possible choice instead of $D_Q$.
Then $D_Q' = D_Q + \ddiv (h)$ for some function $h$ with support disjoint
from $D_P$, and $f_{D_Q'} = f_{D_Q} h^n$, which implies
\[
\frac {f_{D_P} (D_Q')}{f_{D_Q'} (D_P)}
= \frac {f_{D_P} (D_Q) f_{D_P} (\ddiv h)}{f_{D_Q} (D_P) h (D_P)^n}
= \frac {f_{D_P} (D_Q) f_{D_P} (\ddiv h)}{f_{D_Q} (D_P) h (\ddiv f_{D_P})}
= \frac {f_{D_P} (D_Q)}{f_{D_Q} (D_P)}
\]
by Weil reciprocity \eqref {eq:weilrecdisjoint}.
By symmetry, the same argument holds when $D_P$ is chosen differently.

\begin {proofof}{ of equivalence between the second and third definitions}
\label {p:equiv2and3}
A proof is given in \cite[Prop.~8]{Miller04}. The basic idea is to choose
$D_P = [P - R] - [-R]$ and $D_Q = [Q + R] - [R]$ for some point $R$. Then
\eqref {eq:weil3} becomes
\[
\frac {f_{D_P} (Q+R)}{f_{D_Q} (P-R)}
\cdot \frac {f_{D_Q} (-R)}{f_{D_P} (R)}.
\]
Informally, letting $R \to \Oc$, the first factor tends to $\weil (P, Q)$ as
defined in \eqref {eq:weil2}, the second factor tends to $(-1)^n$. This can
be made rigorous using formal groups or the Deuring lift of $E$ to the field
of complex numbers.

Alternatively, one may again use generalised Weil reciprocity. Let $D_P = [P]
- [\Oc]$, so that $f_{D_P} = f_P$. Let $R$ be such that $D_Q = [Q+R] - [R]$
has disjoint support with $D_P$; then $D_Q = [Q] - [\Oc] + \ddiv (h)$ with
$h$ monic at $\Oc$ such that $\ddiv h = [Q + R] - [Q] - [R] + [\Oc]$,
and $f_{D_Q} = f_Q h^n$.

Assume first that $P \neq Q$. Then by Theorem~\ref {th:reciprocity},
\[
1 = \prod_{S \in E (\overline K)} \langle f_P, h \rangle_S
= \frac {f_P (Q + R)}{f_P (R) f_P (Q) h^n (P)}
\cdot (-1)^n
\underbrace {(f_P h^n) (\Oc)}_{= \lc (f_P)}.
\]
So
\begin {eqnarray*}
\frac {f_{D_P} (D_Q)}{f_{D_Q} (D_P)}
& = & \frac {(f_Q h^n) (\Oc)}{(f_Q h^n) (P)}
   \cdot \frac {f_P (Q + R)}{f_P (R)}
= \frac {\lc (f_Q) f_P (Q)}{f_Q (P)}
   \cdot \frac {f_P (Q + R)}{f_P (Q) h^n (P) f_P (R)} \\
& = & (-1)^n \frac {f_P (Q)}{f_Q (P)}
   \cdot \frac {\lc (f_Q)}{\lc (f_P)}
\end {eqnarray*}
by the previous equation.

If $P = Q$, a similar computation shows that \eqref {eq:weil3} evaluates
to~$1$.
\end {proofof}

\section {Tate pairing}
\label {sec:tate}

The Tate pairing has been introduced to cryptology in \cite{FrRu94} as a
means of transporting the discrete logarithm problem from curves into
the multiplicative groups of finite fields.
It goes back to Tate, who in \cite{Tate57} considers abelian varieties
defined over local fields and defines a non-degenerate pairing involving
Galois cohomology groups of the variety and the dual abelian variety.
Lichtenbaum defines in \cite{Lichtenbaum69} a pairing in terms of Picard
groups of curves defined over local fields and their Galois cohomology.
This pairing turns out to be a special case of the Tate pairing and as
such is non-degenerate. Its advantage is that it can easily be computed
in terms of divisors and functions on the curve as stated in~\eqref{eq:tate1}.
See also \cite[\S\S5--8]{Silverman10} for an accessible presentation
of these Galois cohomology related pairings.
By considering torsion elements in the groups and reducing modulo the
discrete valuation of the local field, Frey and Rück obtain a non-degenerate
pairing for curves defined over finite fields.
It is often called the Tate--Lichtenbaum pairing
\cite[\S3.3]{Frey01},\cite[\S6.4.1]{CoFrAvDoLaNgVe06}, although the name
Frey--Rück--Tate--Lichtenbaum pairing might be more appropriate.
In the cryptologic literature, the shorter term Tate pairing has
imposed itself, and we will stick to this tradition.

Computationally, the Tate pairing can be seen as ``half a Weil pairing'';
the idea is to define it directly as $f_P (Q)$ instead of the quotient
\eqref {eq:weil2}. Its precise definition depends on a field extension~$L$
of~$K$ such that $E [n]$ is contained in~$E (L)$; usually, but not
necessarily, $L$ is chosen minimal with this property.

\paragraph {First definition of the Tate pairing.}
Let $P \in E [n]$, let $D_P$ be a degree zero divisor, defined over~$L$,
with $D_P \sim [P] - [\Oc]$, and let $f_{D_P}$, defined over~$L$, be such
that $\ddiv f_{D_P} = n D_P$.
Let $Q$ be another point on $E (L)$ (not necessarily of $n$-torsion) and
let $D_Q \sim [Q] - [\Oc]$ be defined over~$L$ of support disjoint with $D_P$.
Then the Tate pairing of $P$ and $Q$ is given by
\begin {equation}
\label {eq:tate1}
\tate (P, Q) = f_{D_P} (D_Q).
\end {equation}

Algorithm~\ref {alg:miller} shows that $f_{D_P}$ may indeed be defined
over~$L$, so that the pairing takes values in~$L$.
Notice that $f_{D_P}$ is defined only up to a multiplicative constant,
but that this does not change the pairing value since $D_Q$ is of degree~$0$.
Weil reciprocity \eqref {eq:weilrecdisjoint}
shows that if $D_Q$ is replaced by
$D_Q' = D_Q + \ddiv h \sim D_Q$, then \eqref {eq:tate1} is multiplied by
$h (D_P)^n$. Replacing $D_P$ by $D_P' = D_P + \ddiv h$ changes
$f_{D_P}$ to $f_{D_P'} = f_{D_P} h^n$ and thus multiplies the pairing
value by an $n$-th power. So the pairing value is well defined up to
$n$-th powers in~$L$.

Finally, if $Q$ is replaced by $Q + n R$ with $R \in E (L)$, the value
changes again by an $n$-th power. This leads to
adapting the range and domain of $\tate$ as follows.

\begin {theorem}
\label {th:tate}
For $E [n] \subseteq E (L)$, the Tate pairing is a map
\[
\tate : E [n] \times E (L) / n E (L)
\to L^\ast / \left( L^\ast \right)^n
\]
satisfying the following properties
as defined in Theorem~\ref {th:weil}:
\begin {enumerate}
\item
Bilinearity,
\item
Non-degeneracy,
\item
Compatibility with isogenies.
\end {enumerate}
\end {theorem}

\begin {proof}
Bilinearity is immediate from the definition using
$[Q_1 + Q_2] - [\Oc] \sim [Q_1] + [Q_2] - 2 [\Oc]$
by Theorem~\ref {th:principality}, so that
$D_{Q_1 + Q_2} = D_{Q_1} + D_{Q_2}$ and
$f_{P_1 + P_2} = f_{P_1} f_{P_2}$.

Non-degeneracy does not hold over arbitrary fields. In particular, the
pairing becomes completely trivial if every element of $L$ is an $n$-th
power, for instance if $L = \overline K$. So the proofs of non-degeneracy
use the structure of the groups over a finite field, see
\cite{FrRu94,Hess04,Schaefer05,Bruin11}.

Let $\alpha$ be an isogeny. We may assume that $D_P$ and $D_Q$ are chosen
so that all function values encountered during the proof are defined and
non-zero.
From the observation that $D_{\alpha (P)} = \alpha_\ast (D_P)$, one shows
as in~\eqref {eq:divisorisogeny} that
\[
\tate (\alpha (P), \alpha (Q))
= f_{D_{\alpha (P)}} (D_{\alpha (Q)})
= \left( \prod_{R \in \ker \alpha} f_{D_P} \big( (\tau_R)_\ast (D_Q) \big)
\right)^{e_\alpha};
\]
the constant~$c$ of \eqref {eq:divisorisogeny} disappears since $f_{D_P}$
is evaluated in divisors of degree~$0$.
Now Theorem~\ref {th:principality} shows that
$(\tau_R)_\ast (D_Q) \sim D_Q$, so that each factor equals
$\tate (P, Q)$, which finishes the proof.
\end {proof}

Again, an alternative definition yields a computationally advantageous
form of the pairing.

\paragraph {Second definition of the Tate pairing.}
For $P \in E [n]$ and $Q \in E (L)$ (representing a class modulo $n E (L)$),
$P$, $Q \neq \Oc$ and $P \neq Q$, let $f_P$ be monic at $\Oc$ such that
$\ddiv (f_P) = n [P] - n [\Oc]$.
Then
\begin {equation}
\label {eq:tate2}
\tate (P, Q) = \frac {f_P (Q)}{\lc (f_P)};
\end {equation}
if $f_P$ is chosen monic as in Definition~\ref {def:monic},
\[
\tate (P, Q) = f_P (Q).
\]
For one or both of $P$ and $Q$ equal to $\Oc$, one has
$\tate (P, Q) = 1$.
If $P = Q$, one may choose some point $R \in E (L)$ such that
$n R \not\in \{ \Oc, -Q \}$, if it exists, and replace $Q$ by $Q+nR$.

\begin {proofof}{ of equivalence of the two definitions}:
Letting $D_Q = [Q] - [\Oc]$, so that $f_{D_Q} = f_Q$,
and $D_P = [P+R] - [R]$ so that $D_P$ and $D_Q$ have disjoint supports
and $f_{D_P} = f_P h^n$ for the function $h$, monic at $\Oc$, with
$\ddiv (h) = [P+R] - [P] - [R] + [\Oc]$, we immediately obtain
\[
f_{D_P} (D_Q) = \frac {(f_P h^n) (Q)}{(f_P h^n) (\Oc)}
= \frac {f_P (Q) h^n (Q)}{\lc (f_P)}
= \frac {f_P (Q)}{\lc (f_P)}
\]
up to $n$-th powers.
\end {proofof}

Unlike the Weil pairing, the Tate pairing is neither alternating nor
identically~$1$ on the diagonal (which is hardly surprising given that its
two arguments live in different sets). On single $n$-torsion points $P$, it
may or may not hold that $\tate (P, P) = 1$.

The definition of the domain of the Tate pairing as a quotient group is
unwieldy in cryptographic applications, where unique representatives of
pairing results are desired. It can be remedied by observing that $L^\ast$ is
a cyclic group of order $\#L - 1 = q^k - 1$, which is divisible by $n$; so
the map
\[
L^\ast / \left( L^\ast \right)^n \to \mu, \quad
x \mapsto x^{\frac {q^k - 1}{n}}
\]
is an isomorphism with the $n$-th roots of unity $\mu$, and the \textit
{reduced Tate pairing}
\begin {equation}
\label {eq:tatered}
\tatered : E [n] \times E (L) / n E (L) \to \mu, \quad
(P, Q) \mapsto \tate (P, Q)^{\frac {q^k - 1}{n}}
= f_P (Q)^{\frac {q^k-1}{n}}
\end {equation}
(for $P$, $Q \neq \Oc$, $P \neq Q$)
is an equivalent pairing with the same properties as the Tate pairing itself.

It is not generically possible to similarly replace the set $E (L) / n E (L)$
from which the second argument is taken by $E [n]$. As an abelian group, $E
(L)$ is isomorphic to $\Z / r_1 \Z \times \Z / r_2 \Z$ with $n \mid r_1 \mid
r_2$, and $E (L) / n E (L) \simeq \Z / n \Z \times \Z / n \Z$. Consider the
homomorphism
\[
\psi : E (L) / n E (L) \to E [n], \quad
Q \mapsto \frac {r_2}{n} Q.
\]
This homomorphism is injective (and thus an isomorphism by cardinality
considerations) if and only if $\gcd \left( \frac {r_2}{r_1}, n \right) = 1$.
A sufficient (but not necessary) condition is that $\gcd \left( \frac
{r_2}{n}, n \right) = 1$, or equivalently $\gcd \left( \frac {\# E (L)}{n^2},
n \right) = 1$; this is often satisfied in cryptography, where $n$ is a large
prime.
Then the function
\[
e : E [n] \times E [n] \to \mu, \quad
(P, Q) = f_P (Q)^{\frac {q^k-1}{n}}
\]
satisfies $e (P, Q) = \tatered (P, \psi^{-1} (Q))^{\frac {r_2}{n}}$, and
since powering by $\frac {r_2}{n}$ induces a permutation on $\mu$,
it inherits the properties of the reduced Tate pairing.

\section {Computation}
\label {sec:computation}

The main ingredients of the Weil and the Tate pairings are functions with
given divisors; an algorithm computing them is published in \cite{Miller04}
and has become known as Miller's algorithm. The basic idea is to have the
tangent-and-cord law of \S\ref {ssec:divisors} not only reduce a sum of two
points to only one point, but at the same time output the lines that have
served for the reduction. Applied iteratively, it thus reduces a principal
divisor to~$0$ and returns the function having this divisor as a quotient of
products of lines. The algorithm is applicable to any principal divisor, but
we only present it for the case of $n D = n [P] - n [\Oc]$ where $P$ is an
$n$-torsion point, which can be used for computing the Weil pairing
via~\eqref {eq:weil2} and the (reduced) Tate pairing via~\eqref {eq:tate1}
or~\eqref {eq:tate2} and~\eqref {eq:tatered}.

\begin {definition}
\label {def:fi}
For $i \in \Z$, let $f_{i, P}$ be the function (monic at $\Oc$) with
divisor $i [P] - [i P] - (i-1) [\Oc]$.
\end {definition}

The function $f_{i, P}$ exists by Theorem~\ref {th:principality}.
Notice that $f_{1, P} = 1$ and $f_{n, P} = f_P$.
The tangent-and-chord law, applied to $i P$ and $j P$, shows that
\begin {equation}
\label {eq:fsum}
f_{i+j, P} =
f_{i, P} f_{j, P} \frac {\ell_{iP, jP}}{v_{(i+j)P}}
\end {equation}
with $\ell$, $v$ defined as in \eqref {eq:l}, \eqref{eq:v} for $i \not\equiv
-j \pmod n$, $\ell_{i P, (n - i) P} = v_{i P}$ and $v_{\Oc} = 1$. Moreover,
\[
f_{-i, P} = \frac {1}{f_{i, P} v_{iP}}.
\]
These observations yield the following algorithm:

\begin {algorithm}
\label {alg:miller}
\inoutput {An integer $n$ and an $n$-torsion point $P$}{$\ell$ and $v$,
products of lines, such that $f_P = \frac {\ell}{v}$}
\begin {enumerate}
\item
Compute an addition-negation chain $r_1, \ldots, r_s$ for $n$, that is, a
sequence such that $r_1 = 1$, $r_s = n$ and each element $r_i$ is either
\begin {itemize}
\item
the negative of a previsously encountered one: There is $1 \leq j (i) < i$
such that $r_i = - r_{j (i)}$; or
\item
the sum of two previously encountered ones: There are $1 \leq j (i) \leq k
(i) < i$ such that $r_i = r_{j (i)} + r_{k (i)}$.
\end {itemize}
\item
$P_1 \leftarrow P$, $L_1 \leftarrow 1$, $V_1 \leftarrow 1$
\item
\keyword {for} $i = 2, \ldots, s$ \\
\hspace* {5mm} $j \leftarrow j (i), k \leftarrow k (i)$ \\
\hspace* {5mm} \keyword {if} $r_i = - r_j$ \\
\hspace* {10mm} $P_i \leftarrow -P_j$ \\
\hspace* {10mm} $L_i \leftarrow V_j$ \\
\hspace* {10mm} $V_i \leftarrow L_j v_{P_i}$ \\
\hspace* {5mm} \keyword {else} \\
\hspace* {10mm} $P_i \leftarrow P_j + P_k$ \\
\hspace* {10mm} $L_i \leftarrow L_j L_k \ell_{P_{j(i)}, P_{k (i)}}$ \\
\hspace* {10mm} $V_i \leftarrow V_j V_k v_{P_i}$
\item
\keyword {return} $\ell = L_s$, $v = V_s$
\end {enumerate}
\end {algorithm}

Throughout the loop, we have $P_i = r (i) P$ and $\frac {L_i}{V_i} = f_{r
(i), P}$, which proves the correctness of the algorithm. The numerator $\ell$
and the denominator $v$ are computed separately to avoid costly divisions in
a direct computation of $f_P$. Memory handling of the algorithm is simplified
if the standard double-and-add addition chain is used, in which $r_i = 2
r_{i-1}$ or $r_i = r_{i-1} + 1$, so that the result can be accumulated in
two variables $\ell$ and $v$, see \cite[Alg.~IX.1]{Galbraith05}.

For a reasonable addition-negation-chain of length $s \in O (\log n)$, the
algorithm carries out $O (\log n)$ steps. Unfortunately, the degrees of $L_i$
and $V_i$ grow exponentially to reach $O (n)$. This problem can be solved in
two ways: Instead of storing $L_i$ and $V_i$ as dense polynomials, store them
in factored form as a product of lines. This may make sense if several
pairings with the same $P$ are computed.

Otherwise, if $f_P (E)$ is sought for a divisor~$E$, one may compute directly
$L_i (E)$ and $V_i (E)$ during the loop, thus manipulating only elements of
the finite field~$L$; one should then separate again according to the points
with positive or negative multiplicity in $E$ to avoid divisions.
This approach fails when $E$ contains any of the points $P_i = r (i) P$
encountered during the algorithm, which will then be zeroes of some of the
lines. The solution given in \cite {Miller04} is to work with the leading
coefficients of the lines with respect to their Laurent series in local
parameters associated to the points in the support of~$E$
(analogously to Definition~\ref {def:monic}).
Alternatively, one might regroup quotients of consecutive lines having
$P_i$ as zeroes and replace them (by working modulo the curve equation)
by a rational function that is defined and non-zero in $P_i$. Both approaches
are not very practical, since they replace simple arithmetic in a finite
field by more complicated symbolic algebra.
A simpler technique is to replace the divisor~$E$ by an equivalent divisor
not containing any of the $P_i$ in its support, and using~\eqref {eq:weil3}
and~\eqref {eq:tate1}; the price to pay is that~$E$ then contains at least
two points instead of only one in~\eqref {eq:weil2} and~\eqref {eq:tate2}.
Concerning the Tate pairing, since the second argument~$Q$ is defined only
up to $n$-th multiples, one may replace it by $Q + n R$ for some point~$R$.
Finally, one may simply use an addition-negation chain avoiding the support
of~$E$. Since any addition chain necessarily passes through~$2$, it may be
necessary to use negation if $E$ contains~$2 P$ in its support.

The reduced Tate pairing \eqref {eq:tatered} is usually faster to compute
than the Weil pairing \eqref {eq:weil2}: It requires only one instead of two
applications of Algorithm~\ref {alg:miller}. On the other hand, the advantage
is partially lost through the final exponentiation in the reduced Tate
pairing.

\section {Pairings on cyclic subgroups}
\label {sec:subgroups}

All supposedly hard problems on which pairing-based cryptographic primitives
rely can be broken by computing discrete logarithms arbitrarily in $E [n]$ or
the group $\mu$ of $n$-th roots of unity in the embedding field $L$. So
algorithms using Chinese remaindering for discrete logarithms imply that $n =
r$ being prime is the best choice. Then $E [r]$ is a group of order $r^2$
isomorphic to $\Z / r \Z \times \Z / r \Z$. For the sake of security proofs,
it may be desirable to restrict the Weil and reduced Tate pairings to
subgroups, yielding pairings
\[
e : G_1 \times G_2 \to \mu \subseteq L
\]
on cyclic groups $G_i \subset E [r]$ of prime order $r$. In practice, there
is no definite need for such a restriction: The choice of points when
executing the protocol (for instance, by hashing into $E [r]$) implicitly
defines cyclic subgroups $G_i$ generated by these points; but the subgroups
change with each execution of the algorithm. Notice, however, that some
optimised pairings (see \S\ref {sec:shortloop})
can only be defined on specific subgroups, which are
reviewed in the following. An exhaustive description of the cryptographic
properties of different subgroups is given by Galbraith, Paterson and Smart
in \cite{GaPaSm08}. We retain their classification into type~1, 2 and 3
subgroups and pairings and concentrate on the main characteristics
of the different choices.

For the sake of computational efficiency in Algorithm~\ref {alg:miller}, it
is desirable that $G_1$ and $G_2$ be defined over fields that are as small as
possible. So the curve $E (K)$ is chosen such that $r \mid \# E (K)$, and
$G_1$ is generated by a point of order $r$ defined over $K$. As usual in
cryptography, we assume that $k \geq 2$. Then $G_1$ is defined uniquely as
$E (K)[r]$, and the pairing types differ in their selection of $G_2$. An
important cryptographic property that may or may not be given is hashing into
the different groups, or the (essentially equivalent) possibility of random
sampling from the groups. It is a trivial observation that if $H : \{ 0, 1
\}^\ast \to \{ 0, \ldots, r-1 \}$ is a collision-resistant hash-function and
$G_i = \langle P_i \rangle$, then $H_i : \{ 0, 1 \}^\ast \to G_i$,
$m \mapsto H (m) P_i$, is also collision-resistant. But $H_i$ reveals
discrete logarithms,
which breaks most pairing-based cryptographic primitives. A comparatively
expensive way of hashing into $G_1$ is to first hash into a point on $E (K)$
(by hashing to its $X$- or $Y$-coordinate and solving the resulting equation
for the other coordinate; if no solution exists, one needs to hash the
message concatenated with a counter that is increased upon each unsuccessful
trial). One may then multiply by the cofactor $h = \frac {\# E (K)}{r}$,
which yields a point in $G_1$. A similar procedure hashes to arbitrary
$r$-torsion points in~$E (L)$, but these need not lie in a fixed
subgroup~$G_2$.

\subsection {Type 1: \texorpdfstring {$G_1 = G_2$}{G1=G2}}
\label {ssec:type1}

Most of the early papers on pairing-based cryptography are formulated only
for the case of a \textit {symmetric pairing}, in which $G_2 = G_1$. However,
it is in fact not possible to simply choose the arguments of the pairings of
\S\S\ref{sec:weil} and~\ref{sec:tate} from $G_2 = G_1$, since then the
pairing becomes trivial. This is clear for the Weil pairing from
Theorem~\ref {th:weil}(b), but also holds for the reduced Tate pairing:
Algorithm~\ref {alg:miller} implies that the result lies in the field~$K$
over which both
arguments are defined, but $K \cap \mu = \{ 1 \}$. A symmetric pairing may be
obtained for supersingular curves with a so-called \textit {distortion map},
an explicit monomorphism $\psi : E (K)[r] \to E [r] \backslash G_1$. The
non-degeneracy of the Weil pairing then implies that
\[
e : G_1 \times G_1 \to \mu, \quad
(P, Q) \mapsto \rweil (P, \psi (Q))
\]
is also a non-degenerate pairing; the same usually holds for the reduced Tate
pairing.

Algebraic distortion maps cannot exist for ordinary curves, whose
endomorphism rings are abelian. Then $\psi$ would be an
endomorphism and it would commute with the Frobenius, so the image of
$G_1 \subseteq E (K)[r]$ would again lie in $E (K)$ and thus be
equal to $G_1$.

Conversely, supersingular curves have a non-abelian endomorphism ring, and it
has been shown by Galbraith and Rotger in \cite[Th.~5.2]{GaRo04} that they
always admit an algebraic distortion map coming from the theory of complex
multiplication (cf. \cite{Deuring41}) as long as $r \geq 5$; the same article
describes an algorithm for explicitly determining such a map. It is
well-known that supersingular curves with $k = 2$ admit particularly simple
distortion maps, namely,
\begin {equation}
\label {eq:distortion4}
\psi (x, y) = (-x, i y)
\end {equation}
for $E : Y^2 = X^3 + X$ over $\F_p$ with $p \equiv 3 \pmod 4$ and
\begin {equation}
\label {eq:distortion3}
\psi (x, y) = (\zeta_3 x, y)
\end {equation}
for $E : Y^2 = X^3 + 1$ over $\F_p$ with $p \geq 5$ and $p \equiv 2 \pmod 3$,
where $\zeta_3$ and $i$ are primitive third and fourth roots of unity,
respectively, in $\F_{p^2}$.

If the $X$-coordinate of $\psi$ is defined over $K$ (for instance, in \eqref
{eq:distortion4}, but not in \eqref {eq:distortion3}), it is observed in
\cite{BaKiLySc02} that the computation of the reduced Tate pairing
\[
e (P, Q) = \tatered (P, \psi (Q))
= f_P (\psi (Q))^{\frac {q^k - 1}{r}}
\text { by \eqref {eq:tatered}}
\]
can be simplified by omitting denominators. Indeed, notice that if a pure
addition chain (without subtractions) is used, the denominator $v$ returned
by Algorithm~\ref {alg:miller} is a polynomial in $K [X]$ not involving $Y$;
since $X (\psi (Q)) \in K$, the value $v (Q)$ disappears through the final
exponentiation.

The main drawback of type 1 pairings is the lack of flexibility of the
embedding degree~$k$: Since it is limited to supersingular curves,
we have $k \leq 2$ for curves over fields of characteristic at least~$5$,
$k \leq 4$ over fields of characteristic~$2$ and $k \leq 6$ over
fields of characteristic~$3$ by \cite [Theorem~4.1]{Waterhouse69}.

\subsection {Type 2: \texorpdfstring {$G_2 \hookrightarrow G_1$}{G2->G1}}
\label {ssec:type2}

The pairing is of type 2 when there is an efficiently computable monomorphism
$\phi$ from $G_2$ to $G_1$. In some sense, this is the converse of type~1,
where there is a non-trivial monomorphism from $G_1$ into another $r$-torsion
group. This case, however, is essentially the generic one and available in
supersingular and ordinary curves alike. Let $\pi : (x, y) \mapsto (x^q,
y^q)$ be the Frobenius endomorphism related to the field extension $L / K =
\F_{q^k} / \F_q$. Then $K (E)$ is fixed by $\pi$ or, otherwise said, $G_1$
are the $r$-torsion points that are eigenvectors under~$\pi$ with
eigenvalue~$1$. Hasse's theorem then implies that the $r$-torsion of $E$ is
generated by one point $P$ with eigenvalue~$1$ and another point $Q$ with
eigenvalue~$q$. We now consider the \textit {trace} defined as a map on
points by
\[
\Tr : E (L) \to E (K), \quad
R \mapsto \sum_{i = 0}^{k - 1} R^{\pi^i}.
\]
Since the trace of a point is invariant under $\pi$, it is indeed a point
defined over $K$. We have $\Tr (P) = k P \neq \Oc$ in a cryptographic
context, where~$r$ is much bigger than~$k$, and $\Tr (Q) = Q + q Q +
\cdots + q^{k-1} Q = \frac {q^k - 1}{q - 1} Q = \Oc$ since the order~$r$ of
$Q$ divides $q^k - 1$, but not $q - 1$. If $R$ is any $r$-torsion point, then
$R = a P + b Q$, $\Tr (R) = a k P$ and $Q' = k R - \Tr (R) = k b Q \in
\langle Q \rangle$. Unless $R \in \langle P \rangle$, in which case $Q' =
\Oc$, the element $Q'$ is thus a generator of $\langle Q \rangle$, which can
be found efficiently by a randomised algorithm.

Let $R$ be an arbitrary $r$-torsion point that is a pure multiple of neither
$P$ nor $Q$ (which can be checked using the Weil pairing; in practice, a
random $r$-torsion point satisfies this restriction with overwhelming
probability). Let $G_2 = \langle R \rangle$, and $\phi = \Tr$.

The existence of $\phi$ reduces problems (for instance, the discrete
logarithm problem or the decisional Diffie--Hellman problem) defined in terms
of~$G_2$ into problems defined in terms of~$G_1$, which may be helpful for
reductionist security proofs. But as usual, the existence of additional
algebraic structures (here, the map $\phi$) raises doubts as to the
introduction of a security flaw. Furthermore, hashing or random sampling in
$G_2$ appears to be impossible, except for the trivial approach revealing
discrete logarithms. Recent work by Chatterjee and Menezes \cite {ChMe11}
introduces a heuristic construction to transform a cryptographic primitive
in the type 2 setting, together with its security argument, into an
equivalent type 3 primitive. Thus, type 2 pairings should probably be
avoided in practice.

\subsection {Type 3}
\label {ssec:type3}

The remaining case where there is no apparent efficiently computable
mono\-mor\-phism $G_2 \to G_1$ is called type 3. In view of the
discussion of \S\ref {ssec:type2}, this implies that
\begin {eqnarray*}
G_2 & = & \{ R \in E [r] : R^\pi = q R \} \\
    & = & \{ R \in E [r] : \Tr (R) = \Oc \}.
\end {eqnarray*}
The previous discussion has also shown how to find a generator of $G_2$.
Hashing into $G_2$ may be accomplished in a similar manner: Hash to an
arbitrary point $R \in E [r]$, and define $k R - \Tr (R)$ as the final hash
value.

\section {Loop-shortened pairings}
\label {sec:shortloop}

Subsequent work has concentrated on devising pairings with a shorter loop in
Algorithm~\ref {alg:miller}, generally starting from the Tate
pairing~\eqref {eq:tate2}. It turns out that in certain special cases,
\[
e (P, Q) = f_{\lambda, P} (Q)
\text { or }
e (P, Q) = f_{\lambda, Q} (P)
\]
define non-degenerate, bilinear pairings for $\lambda \ll n$ with
$f_{\lambda, P}$ as in Definition~\ref {def:fi}.
The proof proceeds by showing that the pairing is the $M$-th power of the
original Tate pairing for some $M$ prime to~$n$.
Cryptographic applications may then directly use the new pairing, or, for the
sake of interoperability, the Tate pairing may be retrived by an additional
exponentiation with $M^{-1} \bmod n$.
The first such pairing, called $\eta$ pairing, was described by Barreto,
Galbraith, \'O'h\'Eigeartaigh and Scott in \cite {BaGaOhSc07}. It was limited
to supersingular curves and thus yielded a type~1 pairing
(see~\S\ref {ssec:type1}). The examples in \cite {BaGaOhSc07} show that
$\lambda \approx \sqrt n$ is achievable in supersingular curves over fields
of characteristic~$2$ and~$3$.

In the remainder of this section, we fix the same setting as in
\S\ref {sec:subgroups}. In particular, $n = r$ is prime. All pairings will be
defined on $G_1 \times G_2$, where $G_1 = E (K)[r]$ and $G_2$ is the set of
$r$-torsion points defined over $L = \F_{q^k}$ with eigenvalue~$q$ under
the Frobenius $\pi : (x, y) \mapsto (x^q, y^q)$. This is crucial for the
proofs, and incidentally leads to type~3 pairings.

\begin {lemma}
\label {lm:multiplen}
Let $P \in E [n]$.
If $N$ is such that $n \mid N \mid q^k - 1$, then
\[
f_{N, P} = f_{n, P}^{N/n}.
\]
If $N$ is such that $n \mid N$, then
\[
f_{N + 1, P} = f_{N, P}.
\]
\end {lemma}

Both properties hold by definition; the first one was used in
\cite [\S6]{GaHaSo02} to speed up the computation by replacing~$r$ with a small
multiple of low Hamming weight.

\subsection {Ate pairing}

The ate pairing is defined in \cite [Theorem~1]{HeSmVe06} as
\begin {equation}
\label {eq:ate}
\ate : G_1 \times G_2 \to L^\ast / (L^\ast)^r, \quad
(P, Q) \mapsto f_{T, Q} (P)
\end {equation}
with $T = t - 1$, where $t$ is the trace of Frobenius satisfying
$\# E (K) = q + 1 - t$.
 
\begin {theorem}
\label {th:ate}
$\ate$ is bilinear, and if $r^2 \nmid T^k - 1$, it is non-degenerate.
More precisely,
\[
\left( \ate (P, Q) \right)^{k q^{k-1}} = \rtate (Q, P)^{\frac {T^k - 1}{r}}.
\]
\end {theorem}

For the ate pairing and all other pairings presented in the following, a
reduced variant with unique values in $\mu \subseteq L^\ast$ is obtained as
in~\eqref {eq:tatered} by raising to the power $\frac {q^k - 1}{r}$.

\begin {proofof}{ of Theorem~\ref {th:ate}}
The crucial step is the observation that for any $\lambda$,
\begin {eqnarray}
\nonumber
f_{\lambda, T^i Q} \circ \pi^i
& = & f_{\lambda, q^i Q} \circ \pi^i \text { since } T \equiv q \pmod r \\
\nonumber
& = & f_{\lambda, \pi^i (Q)} \circ \pi^i \text { since } Q \in G_2 \\
& = & f_{\lambda, Q}^{q^i},
\label {eq:ffrob}
\end {eqnarray}
since the coefficients of the rational function $f_{\lambda, Q}$ can be
expressed in the coefficients of~$Q$ and of the curve, and the latter lie
in~$\F_q$.

In particular for $P \in G_1$ and $\lambda = T$,
$f_{T, T^i Q} (P) = f_{T, Q}^{q^i} (P)$.

Then
\begin {eqnarray*}
\rtate (Q, P)^{\frac {T^k - 1}{r}}
& = & f_{r, Q}^{\frac {T^k - 1}{r}} (P)
  =   f_{T^k - 1, Q} (P) \text { by Lemma~\ref {lm:multiplen}} \\
& = & f_{T^k, Q} (P) \text { by Lemma~\ref {lm:multiplen} since }
         T^k - 1 \equiv q^k - 1 \equiv 0 \pmod r \\
& = & \prod_{i = 0}^{k - 1} f_{T, T^i Q}^{T^{k - 1 - i}} (P)
         \text { by comparing divisors and collapsing} \\
&& \text {the telescopic sum} \\
& = & f_{T, Q}^{\sum_{i = 0}^{k - 1} T^{k - 1 - i} q^i} (P)
         \text { by \eqref {eq:ffrob}} \\
& = & \ate (P, Q)^{k q^{k - 1}} \text { in } L^\ast / (L^\ast)^r,
         \text { since } T \equiv q \pmod r.
\end {eqnarray*}
\end {proofof}

By Hasse's theorem, $T \in O (\sqrt q)$, so that the number of operations in
Algorithm~\ref {alg:miller} drops generically by a factor of about~$2$; the
effect can, however, be much more noticeable for certain curves.
For instance, \cite {FrScTe10} describes a family of curves for $k = 24$
with $r \in \Theta (q^{4/5})$ and $T \in O (q^{1/10}) = O (r^{1/8})$.
Notice that $8 = \phi (24)$, cf.~\S\ref {ssec:optimal}.
A price to pay is that the arguments~$P$ and~$Q$ are swapped: The elliptic
curve operations need to be carried out over~$\F_{q^k}$ instead of~$\F_q$.
(Algorithm~\ref {alg:miller} in this context is sometimes called
``Miller full'', while the more favourable situation is called
``Miller light''.)

\subsection {Twisted ate pairing}

The twisted variant of the ate pairing keeps the usual order of the
arguments, but sacrifices on the loop length.

Assume $\charac \F_q \geq 5$, and let $d = \gcd (k, \# \Aut (E))$
and $e = \frac {k}{d}$. Then there is a twist~$E'$ of degree~$d$ of~$E$,
that is, a curve~$E'$ defined over~$\F_q$ with an isomorphism
$\psi : E' \to E$, which is defined over $\F_{q^d}$. It can be given
explicitly as follows for $E : Y^2 = X^3 + a X + b$ in short Weierstra{\ss}
form, see \cite [{\S}X.5.4]{Silverman86}:
\[
\begin {array}{lll}
d = 2: &
E' : Y^2 = X^3 + D^2 a X + D^3, &
\psi (x, y) = \left( D x, \sqrt {D^3} y \right); \\
d = 4: &
E' : Y^2 = X^3 + D a X, &
\psi (x, y) = \left( \sqrt D x, \sqrt [4]{D^3} y \right); \\
d \in \{ 3, 6 \}: &
E' : Y^2 = X^3 + D b, &
\psi (x, y) = \left( \sqrt [3]D x, \sqrt D y \right);
\end {array}
\]
where $D$ is a non-square in $\F_q$ for $d \in \{ 2, 4 \}$,
a non-cube and square for $d = 3$,
and a non-cube and non-square for $d = 6$.
The formul{\ae} make sense since
for $d = 4$, we have $b = 0$ and $q \equiv 1 \pmod 4$,
while for $d \in \{ 3, 6 \}$, we have $a = 0$ and $q \equiv 1 \pmod 3$.
Up to isomorphism over $\F_q$, the twist is unique for $d = 2$, and there
are two different ones for $d \in \{ 3, 6 \}$ (such that $g D$ or $g^2 D$,
respectively, is a cube for $g$ a generator of $\F_q^\ast$) and $d = 4$
(such that $g D$ or $g^3 D$, respectively, is a fourth power).
One can then show, see \cite [\S\S4-5]{HeSmVe06}, that besides $E$
itself there is a unique
twist~$E'$ of~$E$, defined over $\F_{q^e}$, such that
$r \mid \# E' (\F_{q^e})$. (This uses that $r^2 \nmid \# E (\F_q)$.)
If $G_2' = E' (\F_{q^e}) [r]$, then $G_2 = \psi (G_2')$.
In particular, the $X$-coordinates of the points in~$G_2$ lie
in~$\F_{q^{k/2}}$ for $d$ even, and the $Y$-coordinates lie
in~$\F_{q^{k/3}}$ for $3 \mid d$.

The twisted ate pairing of \cite [{\S}VI]{HeSmVe06} is defined by
\begin {equation}
\label {eq:twistedate}
\twistedate : G_1 \times G_2 \to L^\ast / (L^\ast)^r, \quad
(P, Q) \mapsto f_{T^e, P} (Q).
\end {equation}

Let $\pi' : (x, y) \mapsto \left( x^q, y^q \right)$ be the Frobenius
of~$E'$, and let the endomorphism~$\alpha$ of~$E$ be defined as
$\alpha = \psi \circ (\pi')^e \circ \psi^{-1}$.
Then $\alpha|_{G_2} = \alpha|_{\psi (G_2')} = \id$,
$\alpha^d|_{G_1} = \id$, and thus $\alpha (G_1) \subseteq G_1$.
Since $\psi$ is an isomorphism and $\deg ((\pi')^e) = q^e$, this implies that
$\alpha|_{G_1}$ is multiplication by~$q^e$.
So~$\alpha$ behaves similarly to the Frobenius endomorphism, but with the
roles of~$G_1$ and~$G_2$ reversed and of degree~$q^e$ instead of~$q$:
$G_2$ is the eigenspace of eiganvalue~$1$, and $G_1$ is the eigenspace of
eigenvalue~$q^e$. The same proof as for Theorem~\ref {th:ate} thus carries
through after replacing~$\pi$ by~$\alpha$, $q$ by~$q^e$, $T$ by~$T^e$
and~$k$ by~$d$.

\begin {theorem}
\label {th:twistedate}
$\twistedate$ is bilinear, and if $r^2 \nmid T^k - 1$, it is non-degenerate.
More precisely,
\[
(\twistedate)^{d q^{e (d-1)}} = (\rtate)^{\frac {T^k - 1}{r}}.
\]
\end {theorem}

Generically. one has $T^e = T^{k/d} \in O \left( q^{k / (2d)} \right)$;
as soon as $k > 2d$, so certainly for $k > 12$, the loop becomes larger than
for the standard Tate pairing, which has the same order of arguments.

\subsection {Optimal pairings}
\label {ssec:optimal}

The discovery of the ate pairing based on a function $f_{\lambda, Q}$,
where $\lambda = T$ is not a multiple of the order of~$Q$ , raised the
question of further possible values for~$\lambda$, and on the possibility
of minimising the loop length~$\log_2 \lambda$.
(Strictly speaking, the loop length in Algorithm~\ref {alg:miller} depends
on the addition-negation chain; $\lfloor \log_2 \lambda \rfloor$ measures
the number of doublings in a standard double-and-add chain.)

For $i = 1, \ldots, k-1$, Zhao, Zhang and Huang define in \cite {ZhZhHu08}
the ate$_i$ pairing by
\begin {equation}
\label {eq:atei}
\atei : G_1 \times G_2 \to L^\ast / (L^\ast)^r, \quad
(P, Q) \mapsto f_{T^i \bmod r, Q} (P).
\end {equation}

For a curve with an automorphism of order $d \mid k$ and
$e = \frac {k}{d}$, a twisted version may be defined for
$i = 1, \ldots, d-1$ as
\[
\twistedatei : G_1 \times G_2 \to L^\ast / (L^\ast)^r, \quad
(P, Q) \mapsto f_{T^{e i} \bmod r, P} (Q).
\]
Their bilinearity and non-degeneracy (if $r^2 \nmid T^{i k'}$, where
$k' = \frac {k}{\gcd (k, i)}$ is the order of~$T^i$ modulo~$r$) is proved
as in Theorems~\ref {th:ate} and~\ref {th:twistedate}, after replacing~$\pi$
by~$\pi^i$ or~$\pi'$ by~$(\pi')^i$, respectively.

In \cite {LeLePa09}, for the first time two such pairings were combined:
If $t_1 = t_0 \lambda_1 + \lambda_0$ and
$f_{t_0, Q}$ and $f_{t_1, Q}$
define powers of the Tate pairing $\rtate (Q, P)$,
then so does
\begin {equation}
\label {eq:rate}
f_{\lambda_1, t_0 Q} f _{\lambda_0, Q}
\frac {\ell_{t_0 \lambda_1 Q, \lambda_0 Q}}{v_{t_1 Q}},
\end {equation}
called the R-ate pairing.
The proof relies on the equation
\begin {equation}
\label {eq:fproduct}
f_{t_0 \lambda_1, Q} = f_{t_0, Q}^{\lambda_1} f_{\lambda_1, t_0 Q},
\end {equation}
which is readily verified by comparing divisors,
so that~\eqref {eq:rate} equals the pairing-defining function
$f_{t_1, Q} f_{t_0, Q}^{\lambda_1}$ by~\eqref {eq:fsum}.
Non-degeneracy holds as soon as the exponent with respect to the Tate
pairing, readily computed from the previous equation, is not divisible
by~$r$. The added loop length in the computation of~\eqref {eq:rate}
is $\log_2 (\lambda_1) + \log_2 (\lambda_0)$. Since the computation
of $f_{\lambda_1, t_0 Q}$ and $f_{\lambda_0, Q}$ by
Algorithm~\ref {alg:miller} finishes with $t_0 \lambda_1 Q$ and
$\lambda_0 Q$, the correction factor is obtained as the quotient of lines
from adding these last two points.
Additionally, $t_0 Q$ needs to be computed (which can be done in parallel
with Algorithm~\ref {alg:miller} for $f_{\lambda_0, Q}$ if an
addition-negation sequence passing through both ~$\lambda_0$ and~$t_0$
is used), and an exponentiation with~$\lambda_1$ is needed, which will
usually be negligeable compared to the final exponentiation for obtaining
reduced pairings.

Several examples of curve families are given in \cite {LeLePa09} with
$t_0$, $t_1$ a power of~$T$ and
$\lambda_0$, $\lambda_1 \in O \left( r^{1 / \phi (k)} \right)$.
That this is no coincidence has been shown by He{\ss} in~\cite {Hess08} and
Vercauteren in~\cite {Vercauteren10}, who defined more general pairing
functions, leading to a notion of optimiality that reaches this quantity
$O \left( r^{1 / \phi (k)} \right)$.

\subsubsection {He{\ss} pairings}
\label {sssec:hess}

\begin {theorem}[\cite {Hess08}, Theorem~1]
Let $t = \sum_{i = 0}^{\deg t} t_i Y^i \in \Z [Y]$ and $y$ a primitive
$k$-th root of unity modulo~$r^2$ such that $r \mid t (y)$. Let
$\fhess$ be the function, monic at~$\Oc$, such that
\begin {equation}
\label {eq:divhess}
\ddiv (\fhess)
= \sum_{i = 0}^{\deg t} t_i \left( [y^i Q] - [\Oc] \right).
\end {equation}
Then the He{\ss} pairing
\begin {equation}
\label {eq:hess}
\hess : G_1 \times G_2 \to L^\ast / (L^\ast)^r, \quad
(P, Q) \mapsto \fhess (P),
\end {equation}
is bilinear and, if $r^2 \nmid t (y)$, non-degenerate.
\end {theorem}

\begin {proof}
Let $t (y) = r L$, and rewrite~\eqref {eq:divhess} as
\[
\ddiv (\fhess) = \sum_{i = 0}^{\deg t} t_i y^i [Q]
- \sum_{i = 0}^{\deg t} t_i \left( y^i [Q] - [y^i Q] \right)
- \left( \sum_{i = 0}^{\deg t} t_i + 1 \right) [\Oc],
\]
which implies that
\[
\fhess = f_{r, Q}^L \prod_{i = 0}^{\deg t} \left( f_{y^i, Q} \right)^{- t_i}.
\]
Since~$q$ is a primitive $k$-th root of unity modulo~$r$, we have
$y \equiv q^j \pmod r$ for some~$j$, and $y^i \equiv q^{ij} \pmod r$.
The same proof as for the ate (or ate$_i$) pairing, with $y^i$ in the place
of~$T$ and~$\pi^{ij}$ in the place of~$\pi$, shows that
\[
f_{y^i, Q}^{k q^{k-1}} (P) = \rtate (Q, P)^{\frac {y^{ik} - 1}{r}} = 1
\text { since } r^2 \mid y^k - 1.
\]
Since $r \nmid k q^{k-1}$, we have $f_{y^i, Q} (P) = 1$.
So $\hess = (\rtate)^L$ is bilinear, and non-degenerate for $r \nmid L$.
\end {proof}

\begin {remark}
\label {rem:hess}
The condition that $y$ be a primitive $k$-th root of unity modulo~$r^2$ is
clearly not necessary. If $y$ is a root of unity modulo~$r$, then the
previous proof carries through, showing that $\hess$ is bilinear. More
precisely, $(\hess)^{k q^{k-1}} = (\rtate)^N$ with
\[
N = k q^{k-1} \frac {t (y)}{r}
- \sum_{i = 0}^{\deg t} t_i \frac {y^{ik} - 1}{r}
= \frac {1}{r} \left( k q^{k-1} t (y) - (t (y^k) - t (1)) \right),
\]
so that $\hess$ is non-degenerate if and only if
$r \nmid k q^{k-1} t (y) - \left( t (y^k) - t (1) \right)$.
This should hold with overwhelming probability.
For instance, one can usually choose $y = T = q \bmod r$.
\end {remark}

Since~$y$ is a $k$-th root of unity modulo the order~$r$ of~$Q$, any
function as in~\eqref {eq:divhess} is realised by a polynomial~$t$ of
degree at most $\phi (k) - 1$. Those with a root in~$y$ modulo~$r$ can be
seen as elements of the $\Z$-lattice with basis $r, Y - y,
Y^2 - (y^2 \bmod r), \ldots, Y^{\phi (k) - 1} - (y^{\phi (k) - 1} \bmod r)$
of dimension $\phi (k)$ and determinant~$r$. For fixed dimension, the LLL
algorithm finds an element~$t$ of degree at most $\phi (k) - 1$ and with
$|t_i| \in O \left( r^{1 / \phi (k)} \right)$.

There is a twisted variant of the He{\ss} pairing:
If~$E$ has a twist of order~$d \mid k$ and $e = \frac {k}{d}$, $y$ is a
$d$-th root of unity modulo~$r$ and $r \mid t (y)$, then
\[
\twistedhess : G_1 \times G_2 \to L^\ast / (L^\ast)^r, \quad
(P, Q) \mapsto \ftwistedhess (Q)
\]
defines a bilinear pairing that is non-degenerate if $y$ is a primitive
$d$-th root of unity modulo~$r^2$ or, more generally, if
$r^2 \nmid d q^{e (d-1)} t (y) - \left( t (y^d) - t (1) \right)$.
Using LLL, one obtains a polynomial of degree less than~$\phi (d)$ and
with $|t_i| \in O \left( r^{1 / \phi (d)} \right)$. The only cases of
interest are $d \in \{ 3, 4, 6 \}$, for which $\phi (d) = 2$.
Even then, there is only a constant gain in the loop length that does not
increase with~$k$, so that asymptotically, the He{\ss} pairing will be
preferred to its twisted version.
Finally, \cite {Hess08} also contains an optimal version of the Weil pairing.

To see whether~\eqref {eq:hess} can be computed efficiently, let
$R_i = y^i Q$, $s_i = \sum_{j=0}^i t_j y^j$ and
$S_i = s_i Q = \sum_{j=0}^i t_j R^j$ for $i \geq 0$ and
$s_{-1} = 0$ and $S_{-1} = \Oc$.
Then~\eqref {eq:hess} can be rewritten as
\begin {eqnarray*}
&& \sum_{i=0}^{\deg t} t_i \left( [R_i] - [\Oc] \right) \\
& = & \sum_{i=0}^{\deg t} \ddiv (f_{t_i, R_i})
+ \sum_{i=0}^{\deg t} \left( [t_i R_i] - [\Oc] \right) \\
& = & \sum_{i=0}^{\deg t} \ddiv (f_{t_i, R_i})
+ \sum_{i=0}^{\deg t} \left( [S_i] - [S_{i-1}]
   + \ddiv \left( \frac {\ell_{S_{i-1}, t_i R_i}}{v_{S_i}} \right) \right)
\end {eqnarray*}
and
\[
\fhess =
\prod_{i=0}^{\deg t} f_{t_i, R_i}
\prod_{i=0}^{\deg t} \frac {\ell_{S_{i-1}, t_i R_i}}{v_{S_i}}.
\]
The precomputation of the $R_i$ by $\deg t - 1$ scalar multiplications
can already be rather costly. As $t_i R_i$ is a sideproduct of the
computation of $f_{t_i, R_i}$, each quotient of two lines comes out of a
point addition on $E (L)$.
But by computing each $f_{t_i, R_i}$ separately via
Algorithm~\ref {alg:miller}, the factor $\phi (k)$ gained in the loop length
is lost again through the number of evaluations. So while it is shown in
\cite [Lemma~1]{Hess08} that the He{\ss} pairing uses a function of
relatively low degree in $O \left( r^{1 / \phi (k)} \right)$, it is unclear
whether this function can always be evaluated in
$\frac {\log_2 (r)}{\phi (k)}$ steps or a very small multiple thereof.

\subsubsection {Vercauteren pairings}
\label {sssec:vercauteren}

If one removes the condition that $y$ be a primitive $k$-th root of unity
modulo~$r^2$ in the He{\ss} pairing, one may let $y = q$ under the
conditions of Remark~\ref {rem:hess}, a special case considered
independently by Vercauteren in \cite {Vercauteren10}.
Then the $R_i$ may be computed by successive applications of the Frobenius
map, and moreover,
\[
f_{t_i, R_i} (P)
= f_{t_i, q^i Q} (P)
= f_{t_i, Q}^{q^i} (P) \text { by \eqref {eq:ffrob}}.
\]
These functions have the advantage of being computed by
Algorithm~\ref {alg:miller} with respect to the same base point~$Q$.
By choosing an addition-negation sequence that passes through all
the~$t_i$, they may thus be obtained at the same time.
Currently known algorithms compute such sequences with
$\log_2 N + \phi (k) O \left( \frac {\log N}{\log \log N} \right)$
steps, where $N = \max |t_i|$, for instance by \cite {Yao76}.
This shows that, up to the minor factor $\log \log N$, again the gain of
$\phi (k)$ in the loop lengths is offset by the number of functions.
One should notice, however, that better addition sequences can often be
found in practice. Moreover, coefficients occurring in a pairing context
are far from random, but exhibit arithmetic peculiarities, as illustrated
in the next paragraph.

\subsubsection {Optimal pairings on curve families}

Elliptic curves suitable for pairing-based cryptography, that is, with a
small embedding degree~$k$, are extremely rare among all elliptic curves,
see \cite {Boxall12}. An excellent survey article on the problem of
finding good parameter combinations is \cite {FrScTe10}, so there is no
need to give any details here. Starting with the article by Brezing and
Weng \cite {BrWe05}, work has concentrated on finding families of curves
parameterised by polynomials. For fixed~$k$, these are given by $p (X)$,
$r (X)$ and $u (X) \in \Z [X]$ satisfying arithmetic properties so that
if $x_0 \in \Z$ such that $p (x_0)$ is prime, then there is an elliptic
curve over $\F_{p (x_0)}$ with trace of Frobenius $u (x_0)$ and a
subgroup of order $r (x_0)$ of embedding degree~$k$. Concrete instances
are thus given whenever $p (X)$ and $r (X)$ simultaneously represent
primes. In practice, one has $\deg (p (X)) = \phi (k)$ or $2 \phi (k)$,
and the polynomials tend to have small and arithmetically meaningful
coefficients (for instance, they are often divisible by prime factors
of~$k$).

As an example, Freeman gives a family in \cite [Theorem~3.1]{Freeman06}
for $k = 10$ with
\begin {eqnarray*}
p (X) & = & 25 X^4 + 25 X^3 + 25 X^2 + 10 X + 3, \\
u (X) & = & 10 X^2 + 5 X + 3, \\
r (X) & = & 25 X^4 + 25 X^3 + 15 X^2 + 5 X + 1.
\end {eqnarray*}\
To construct optimal pairings, one may now work directly with polynoials
instead of integers, looking for short vectors in the $\Z [X]$-lattice
with basis
\[
r (X), Y - y (X), Y^2 - \left( y (X)^2 \bmod r (X) \right),
\ldots, Y^{\phi (k)} - \left( y (X)^{\phi (k)} \bmod r (X) \right).
\]

In He{\ss}'s construction of \S\ref {sssec:hess}, $y (X)$ is hereby a
primitive $k$-th root of unity modulo $r (X)^2$; notice that $r (X)$
is necessarily irreducible since it represents primes.

For Vercauteren's specialisation of \S\ref {sssec:vercauteren}, one has
$y (X) = p (X)$, and the above family leads to a short vector
(see \cite [{\S}IV.B]{Vercauteren10})
\[
t (Y) = X Y^3 + X Y^2 - X Y - (X+1).
\]
This means that whenever $p (x_0)$ and $r (x_0)$ are prime for some
$x_0 \in \Z$, then we obtain a curve and an optimal pairing in which the
computation of the $f_{t_i (x_0), Q}$ boils down to $f_{x_0, Q}$. Notice
that $x_0 \approx r (x_0)^{1 / \deg r (X)} = r (x_0)^{1 / \phi (10)}$,
and in this family, the gain of a factor of $\phi (k)$ in each
invocation of Algorithm~\ref {alg:miller} leads indeed to a corresponding
speed-up in the complete function evaluation.

\paragraph {Acknowledgement.}

This research was partially funded by ERC Starting Grant ANTICS 278537.

\bibliographystyle {alpha}
\bibliography {pairings}

\newcommand{\etalchar}[1]{$^{#1}$}
\begin{thebibliography}{TOOO07}

\bibitem[AL13]{AbLa13}
Michel Abdalla and Tanja Lange, editors.
\newblock {\em Pairing-Based Cryptography --- Pairing 2012}, volume 7708 of
  {\em Lecture Notes in Computer Science}, Heidelberg, 2013. Springer-Verlag.

\bibitem[BF01]{BoFr01}
Dan Boneh and Matt Franklin.
\newblock Identity-based encryption from the {W}eil pairing.
\newblock In Joe Kilian, editor, {\em Advances in Cryptology --- CRYPTO 2001},
  volume 2139 of {\em Lecture Notes in Computer Science}, pages 213--229,
  Berlin, 2001. Springer-Verlag.

\bibitem[BGOS07]{BaGaOhSc07}
Paulo S. L.~M. Barreto, Steven~D. Galbraith, Colm \'O'h\'Eigeartaigh, and
  Michael Scott.
\newblock Efficient pairing computation on supersingular abelian varieties.
\newblock {\em Designs, Codes and Cryptography}, 42:239--271, 2007.

\bibitem[BK98]{BaKo98}
R.~Balasubramanian and N.~Koblitz.
\newblock The improbability that an elliptic curve has subexponential discrete
  log problem under the {M}enezes--{O}kamoto--{V}anstone algorithm.
\newblock {\em Journal of Cryptology}, 11:141--145, 1998.

\bibitem[BKLS02]{BaKiLySc02}
Paulo S. L.~M. Barreto, Hae~Y. Kim, Ben Lynn, and Michael Scott.
\newblock Efficient algorithms for pairing-based cryptosystems.
\newblock In Moti Yung, editor, {\em Advances in Cryptology --- CRYPTO 2002},
  volume 2442 of {\em Lecture Notes in Computer Science}, pages 354--369,
  Berlin, 2002. Springer-Verlag.

\bibitem[Box12]{Boxall12}
John Boxall.
\newblock Heuristics on pairing-friendly elliptic curves.
\newblock {\em Journal of Mathematical Cryptology}, 6(2):81--104, 2012.

\bibitem[Bru11]{Bruin11}
Peter Bruin.
\newblock The {T}ate pairing for abelian varieties over finite fields.
\newblock {\em Journal de Th\'eorie des Nombres de Bordeaux}, 23(2):323--328,
  2011.

\bibitem[BSS99]{BlSeSm99}
Ian Blake, Gadiel Seroussi, and Nigel Smart.
\newblock {\em Elliptic Curves in Cryptography}, volume 265 of {\em London
  Mathematical Society Lecture Note Series}.
\newblock Cambridge University Press, Cambridge, 1999.

\bibitem[BW05]{BrWe05}
Friederike Brezing and Annegret Weng.
\newblock Elliptic curves suitable for pairing based cryptography.
\newblock {\em Designs, Codes and Cryptography}, 37(1):133--141, 2005.

\bibitem[CC90]{ChCo90}
Leonard~S. Charlap and Raymond Coley.
\newblock An elementary introduction to elliptic curves {II}.
\newblock CCR Expository Report~34, Institute for Defense Analyses, Princeton,
  July 1990.
\newblock \url {http://www.idaccr.org/reports/er34.ps}.

\bibitem[CFA{\etalchar{+}}06]{CoFrAvDoLaNgVe06}
Henri Cohen, Gerhard Frey, Roberto Avanzi, Christophe Doche, Tanja Lange, Kim
  Nguyen, and Frederik Vercauteren.
\newblock {\em Handbook of Elliptic and Hyperelliptic Curve Cryptography}.
\newblock Discrete mathematics and its applications. Chapman \& Hall/CRC, Boca
  Raton, 2006.

\bibitem[CM11]{ChMe11}
Sanjit Chatterjee and Alfred Menezes.
\newblock On cryptographic protocols employing asymmetric pairings -- the role
  of $\psi$ revisited.
\newblock {\em Discrete Applied Mathematics}, 159:1311--1322, 2011.

\bibitem[CZ14]{CaZh14}
Zhenfu Cao and Fangguo Zhang, editors.
\newblock {\em Pairing-Based Cryptography --- Pairing 2014}, volume 8365 of
  {\em Lecture Notes in Computer Science}, Cham, 2014. Springer-Verlag.

\bibitem[Deu41]{Deuring41}
Max Deuring.
\newblock Die {T}ypen der {M}ultiplikatorenringe elliptischer
  {F}unktionenk\"orper.
\newblock {\em Abhandlungen aus dem mathematischen Se\-mi\-nar der
  ham\-bur\-gi\-schen Universit\"at}, 14:197--272, 1941.

\bibitem[Eng99]{Enge99}
Andreas Enge.
\newblock {\em Elliptic Curves and Their Applications to Cryp\-to\-graphy ---
  {A}n Introduction}.
\newblock Kluwer Academic Publishers, 1999.

\bibitem[Eng08]{Enge08}
Andreas Enge.
\newblock Discrete logarithms in curves over finite fields.
\newblock In Gary~L. Mullen, Daniel Panario, and Igor~E. Shparlinski, editors,
  {\em Finite Fields and Applications}, volume 461 of {\em Contemporary
  Mathematics}, pages 119--139. American Mathematical Society, 2008.

\bibitem[FR94]{FrRu94}
Gerhard Frey and Hans-Georg R\"uck.
\newblock A remark concerning $m$-divisibility and the discrete logarithm in
  the divisor class group of curves.
\newblock {\em Mathematics of Computation}, 62(206):865--874, April 1994.

\bibitem[Fre01]{Frey01}
Gerhard Frey.
\newblock Applications of arithmetical geometry to cryptographic constructions.
\newblock In Dieter Jungnickel and Harald Niederreiter, editors, {\em Finite
  Fields and Applications --- Proceedings of The Fifth International Conference
  on Finite Fields and Applications F$_{q^5}$, held at the University of
  Augsburg, Germany, August 2--6, 1999}, pages 128--161, Berlin, 2001.
  Springer-Verlag.

\bibitem[Fre06]{Freeman06}
David Freeman.
\newblock Constructing pairing-friendly elliptic curves with embedding
  degree~10.
\newblock In Florian Hess, Sebastian Pauli, and Michael Pohst, editors, {\em
  Algorithmic Number Theory --- ANTS-VII}, volume 4076 of {\em Lecture Notes in
  Computer Science}, pages 452--465, Berlin, 2006. Springer-Verlag.

\bibitem[FST10]{FrScTe10}
David Freemann, Michael Scott, and Edlyn Teske.
\newblock A taxonomy of pairing-friendly elliptic curves.
\newblock {\em Journal of Cryptology}, 23(2):224--280, 2010.

\bibitem[Gal05]{Galbraith05}
Steven Galbraith.
\newblock Pairings.
\newblock In Ian~F. Blake, Gadiel Seroussi, and Nigel~P. Smart, editors, {\em
  Advances in Elliptic Curve Cryptography}, chapter~9, pages 183--213.
  Cambridge University Press, Cambridge, 2005.

\bibitem[GHS02]{GaHaSo02}
Steven~D. Galbraith, Keith Harrison, and David Soldera.
\newblock Implementing the {T}ate pairing.
\newblock In Claus Fieker and David~R. Kohel, editors, {\em Algorithmic Number
  Theory --- ANTS-V}, volume 2369 of {\em Lecture Notes in Computer Science},
  pages 324--337, Berlin, 2002. Springer-Verlag.

\bibitem[GP08]{GaPa08}
Steven~D. Galbraith and Kenneth~G. Paterson, editors.
\newblock {\em Pairing-Based Cryptography --- Pairing 2008}, volume 5209 of
  {\em Lecture Notes in Computer Science}, Berlin, 2008. Springer-Verlag.

\bibitem[GPS08]{GaPaSm08}
Steven~D. Galbraith, Kenneth~G. Paterson, and Nigel~P. Smart.
\newblock Pairings for cryptographers.
\newblock {\em Discrete Applied Mathematics}, 156(16):3113--3121, 2008.

\bibitem[GR04]{GaRo04}
Steven~D. Galbraith and Victor Rotger.
\newblock Easy decision {D}iffie--{H}ellman groups.
\newblock {\em LMS Journal of Computation and Mathematics}, 7:201--218, 2004.

\bibitem[Hes08]{Hess08}
Florian Hess.
\newblock Pairing lattices.
\newblock In S.~D. Galbraith and K.~Paterson, editors, {\em Pairing-Based
  Cryptography --- Pairing 2008}, volume 5209 of {\em Lecture Notes in Computer
  Science}, pages 18--38, Berlin, 2008. Springer-Verlag.

\bibitem[He{\ss}04]{Hess04}
Florian He{\ss}.
\newblock A note on the {T}ate pairing of curves over finite fields.
\newblock {\em Archiv der Mathematik}, 82:28--32, 2004.

\bibitem[HSV06]{HeSmVe06}
Florian Hess, Nigel~P. Smart, and Frederik Vercauteren.
\newblock The eta pairing revisited.
\newblock {\em IEEE Transactions on Information Theory}, 52(10):4595--4602,
  2006.

\bibitem[JMO10]{JoMiOt10}
Marc Joye, Atsuko Miyaji, and Akira Otsuka, editors.
\newblock {\em Pairing-Based Cryptography --- Pairing 2010}, volume 6487 of
  {\em Lecture Notes in Computer Science}, Berlin, 2010. Springer-Verlag.

\bibitem[Jou00]{Joux00}
Antoine Joux.
\newblock A one round protocol for tripartite {D}iffie--{H}ellman.
\newblock In Wieb Bosma, editor, {\em Algorithmic Number Theory --- ANTS-IV},
  volume 1838 of {\em Lecture Notes in Computer Science}, pages 385--393,
  Berlin, 2000. Springer-Verlag.

\bibitem[Lic69]{Lichtenbaum69}
Stephen Lichtenbaum.
\newblock Duality theorems for curves over $p$-adic fields.
\newblock {\em Inventiones mathematicae}, 7(2):120--136, 1969.

\bibitem[LLP09]{LeLePa09}
Eunjeong Lee, Hyang-Sook Lee, and Cheol-Min Park.
\newblock Efficient and generalized pairing computation on abelian varieties.
\newblock {\em IEEE Transactions on Information Theory}, 55(4):1793--1803,
  2009.

\bibitem[Mil04]{Miller04}
Victor~S. Miller.
\newblock The {W}eil pairing, and its efficient calculation.
\newblock {\em Journal of Cryptology}, 17:235--261, 2004.

\bibitem[MOV93]{MeOkVa93}
Alfred~J. Menezes, Tatsuaki Okamoto, and Scott~A. Vanstone.
\newblock Reducing elliptic curve logarithms to logarithms in a finite field.
\newblock {\em IEEE Transactions on Information Theory}, 39(5):1639--1646,
  September 1993.

\bibitem[Odl13]{Odlyzko13}
Andrew Odlyzko.
\newblock Discrete logarithms over finite fields.
\newblock In Gary~L. Mullen and Daniel Panario, editors, {\em Handbook of
  Finite Fields}, Discrete Mathematics and Its Applications, chapter 11.6, page
  393–401. Chapman and Hall/CRC, Boca Raton, 2013.

\bibitem[Sch05]{Schaefer05}
Edward~F. Schaefer.
\newblock A new proof for the non-degeneracy of the {F}rey-{R}\"uck pairing and
  a connection to isogenies over the base field.
\newblock In Tanush Shaska, editor, {\em Computational aspects of algebraic
  curves}, volume~13 of {\em Lecture Notes Series on Computing}, pages 1--12,
  Singapore, 2005. World Scientific Publishing Company.

\bibitem[Sil86]{Silverman86}
Joseph~H. Silverman.
\newblock {\em The Arithmetic of Elliptic Curves}.
\newblock Graduate Texts in Mathematics. Springer-Verlag, New York, 2nd 2009
  edition, 1986.

\bibitem[Sil10]{Silverman10}
Joseph~H. Silverman.
\newblock A survey of local and global pairings on elliptic curves and abelian
  varieties.
\newblock In Marc Joye, Atsuko Miyaji, and Akira Otsuka, editors, {\em
  Pairing-Based Cryptography --- Pairing 2010}, volume 6487 of {\em Lecture
  Notes in Computer Science}, pages 377--396, Berlin, 2010. Springer-Verlag.

\bibitem[Sil13]{Silverman13}
Joseph Silverman.
\newblock Elliptic curves.
\newblock In Gary~L. Mullen and Daniel Panario, editors, {\em Handbook of
  Finite Fields}, Discrete Mathematics and Its Applications, chapter 12.2.
  Chapman and Hall/CRC, Boca Raton, 2013.

\bibitem[SOK00]{SaOhKa00}
R.~Sakai, K.~Ohgishi, and M.~Kasahara.
\newblock Cryptosystems based on pairing, 2000.
\newblock SCIS 2000, The 2000 Symposium on Cryptography and Information
  Security, Okinawa, Japan, January 26--28.

\bibitem[SW09]{ShWa09}
Hovav Shacham and Brent Waters, editors.
\newblock {\em Pairing-Based Cryptography --- Pairing 2009}, volume 5671 of
  {\em Lecture Notes in Computer Science}, Berlin, 2009. Springer-Verlag.

\bibitem[Tat58]{Tate57}
J.~Tate.
\newblock {WC}-groups over $p$-adic fields. {E}xposé no. 156.
\newblock In {\em Années 1956/57--1957/58, exposés 137--168}, volume~4 of
  {\em Séminaire Bourbaki}, pages 265--277. Société Mathématique de France,
  1956--1958.

\bibitem[TOOO07]{TaOkOkOk07}
Tsuyoshi Takagi, Tatsuaki Okamoto, Eiji Okamoto, and Takeshi Okamoto, editors.
\newblock {\em Paring-Based Cryptography --- Pairing 2007}, volume 4575 of {\em
  Lecture Notes in Computer Science}, Berlin, 2007. Springer-Verlag.

\bibitem[Ver10]{Vercauteren10}
Frederik Vercauteren.
\newblock Optimal pairings.
\newblock {\em IEEE Transactions on Information Theory}, 56(1):455--461, 2010.

\bibitem[Wat69]{Waterhouse69}
William~C. Waterhouse.
\newblock Abelian varieties over finite fields.
\newblock {\em Annales Scientifiques de l'\'Ecole Normale Sup\'erieure, $4^e$
  S\'erie}, 2:521--560, 1969.

\bibitem[Wei40]{Weil40}
Andr{\'e} Weil.
\newblock Sur les fonctions alg\'ebriques \`a corps de constantes fini.
\newblock {\em Comptes rendus hebdomadaires des s\'eances de l'Aca\-d\'e\-mie
  des sciences}, 210:592--594, 1940.

\bibitem[Yao76]{Yao76}
Andrew Chi-Chih Yao.
\newblock On the evaluation of powers.
\newblock {\em SIAM Journal on Computing}, 5(1):100--103, March 1976.

\bibitem[ZZH08]{ZhZhHu08}
Chang-An Zhao, Fangguo Zhang, and Jiwu Huang.
\newblock A note on the {A}te pairing.
\newblock {\em International Journal of Information Security}, 7(6):379--382,
  2008.

\end{thebibliography}

\end {document}